\documentstyle[leqno]{article}
\newcommand{\se}[1]{{\section{#1}} {\setcounter{equation}{0}}}
\newtheorem{th}{Theorem}[section]
\newtheorem{lm}{Lemma}[section]
\newtheorem{prop}{Proposition}[section]
\newtheorem{de}{Definition}[section]
\newtheorem{co}{Corollary}[section]

\def\k{{K\"{a}hler }}

\def\cy{{Calabi-Yau }}
\def\l{{Lagrangian }}

\begin{document}
\hbadness=10000
\title{{\bf Lagrangian torus fibration and mirror symmetry of Calabi-Yau hypersurface in toric variety}}
\author{Wei-Dong Ruan\\
Department of mathematics\\
Columbia University\\
New York, NY 10027\\}
\footnotetext{Partially supported by NSF Grant DMS-9703870.}
\maketitle
\tableofcontents
\se{Introduction}
In this paper we give a construction of Lagrangian torus fibration for Calabi-Yau hypersurface in toric variety via the method of gradient flow. Using our construction of Lagrangian torus fibration, we are able to prove the symplectic topological version of SYZ mirror conjecture for generic Calabi-Yau hypersurface in toric variety. We will also be able to give precise formulation of SYZ mirror conjecture in general (including singular locus and duality of singular fibres).\\\\
The motivation of our work comes from the study of Mirror Symmetry. Mirror Symmetry conjecture originated from physicists' work in conformal field theory and string theory. It proposes that for a Calabi-Yau 3-fold $X$ there exists a Calabi-Yau 3-fold $Y$ as its mirror manifold. The quantum geometry of $X$ and $Y$ are closely related. In particular one can compute the number of rational curves in $X$ by solving the Picard-Fuchs equation coming from variation of Hodge structure of $Y$.\\\\
Despite the great impact and success mirror symmetry conjecture brings to the understanding of geometry of Calabi-Yau manifolds and their moduli spaces, the fundamental question that in general given a Calabi-Yau manifold how to construct the corresponding mirror Calabi-Yau manifold is not clear at all from the original mirror conjecture. Although mirror manifolds are worked out in certain cases, the construction of mirror for general Calabi-Yau is still very elusive and mysterious. The most general construction so far was given by Batyrev for Calabi-Yau hypersurface in toric variety. From the toric geometry standing point, Batyrev propose that for Calabi-Yau hypersurface $X$ in the toric variety $P_\Delta$ corresponding to a reflexive polyhedron $\Delta$, the mirror should be the Calabi-Yau hypersurface $Y$ in the toric variety $P_{\Delta^\vee}$ corresponding to the dual reflexive polyhedron $\Delta^\vee$. Batyrev computed among other things the Hodge numbers of Calabi-Yau hypersurface $X$ and the mirror $Y$ and showed that they behave as predicted by mirror symmetry. Batyrev's mirror construction includes many special cases of mirror construction discussed previously by many physicists and mathematicians. Later, there are lots of important work surrounding Batyrev mirror construction, including computing number of rational curves, etc.\\\\ 
In 1996 Strominger, Yau and Zaslow (\cite{SYZ}) proposed a geometric construction of mirror manifold via special Lagrangian torus fibration. According to their program (we will call it SYZ construction), a Calabi-Yau 3-fold should admit a special Lagrangian torus fibration. The mirror manifold can be obtained by dualizing the fibers. Or equivalently, the mirror manifold of $M$ is the moduli space of special Lagrangian 3-torus in $M$ with a flat $U(1)$ connection. \\\\
In a sense, SYZ construction appear to be more fundamental to mirror symmetry phenomenon and more classical than quantum mirror symmetry. More importantly, SYZ construction has the potential to explain the mathematical reason behind mirror symmetry. For example, SYZ construction gives us a possible way to construct the mirror manifold if we understand how to construct dual singular fibers.\\\\
{\bf Remark:} The original SYZ mirror conjecture is rather sketchy in nature. Many detail knowledge of the fibraton like singular locus, singular fibres and duality of singular fibres are necessary to be worked out to formulate the {\bf precise SYZ mirror conjecture}. Without the precise formulation, one would not really be able to construct the mirror manifold completely.\\\\
According to the SYZ construction, special Lagrangian submanifold and special Lagrangian fibration for Calabi-Yau manifolds seem to play very important roles in understanding mirror symmetry. However, despite its great potential in solving the mirror symmetry conjecture, our understanding on special Lagrangian submanifolds is very limited. The known examples are mostly explicit local examples or examples coming from $n=2$. There are very few example of special Lagrangian submanifold or special Lagrangian fibration for dimension higher than two. M. Gross, P.M.H. Wilson and N. Hitchin (\cite{Gross1}\cite{Gross2}\cite{GW}\cite{H}) did some important work in this area in recent years. On the other extreme, in \cite{Z}, Zharkov constructed some non-Lagrangian fibration of Calabi-Yau hypersurface in toric variety. Despite all these efforts, SYZ construction still remains to be a beautiful dream to us.\\\\ 
For our discussion, we will relax the special Lagrangian condition, yet still keep the Lagrangian condition. We will consider Lagrangian torus fibration of Calabi-Yau manifolds. According to our discussion in \cite{lag2}, There are major differences between $C^\infty$-Lagrangian fibration and general Lagrangian fibration, let alone non-Lagrangian fibration. In a sense, $C^\infty$-Lagrangian fibration should already capture symplectic topological structure of the corresponding special Lagrangian fibration. We will construct Lagrangian torus fibration that exhibit the same topological structure as $C^\infty$-Lagrangian torus fibration. In particular, singular locus of the fibration will be of codimension 2. In \cite{smooth} we smooth out our Lagrangian fibration to a great extent, yet we still fall short of making it into a $C^\infty$-Lagrangian fibration everywhere.\\\\ 
In this paper, we will construct Lagrangian torus fibration of generic Calabi-Yau hypersurface in toric variety correponding to a reflexive polyhedron in complete generality. With these detailed understanding of Lagrangian torus fibration of generic Calabi-Yau hypersurface in toric variety, we will be able to prove the symplectic topological SYZ mirror conjecture for Calabi-Yau hypersurface in toric variety. More precisely\\
\begin{th}
For generic Calabi-Yau hypersurface $X$ and its mirror Calabi-Yau hypersurface $Y$ near their corresponding large complex limit and large radius limit, there exist corresponding Lagrangian torus fibrations\\
\[
\begin{array}{ccccccc}
X_{\phi(b)}&\hookrightarrow& X& \ \ &Y_b&\hookrightarrow& Y\\
&&\downarrow& \ &&&\downarrow\\
&& \partial \Delta_v& \ &&& \partial \Delta_w\\
\end{array}
\]\\
with singular locus $\Gamma \subset \partial \Delta_v$ and $\Gamma' \subset \partial \Delta_w$, where $\phi:\partial \Delta_w \rightarrow \partial \Delta_v$ is a natural homeomorphism, $\phi(\Gamma')=\Gamma$, and $b\in \partial \Delta_w \backslash \Gamma'$. The corresponding fibres $X_{\phi(b)}$ and $Y_b$ are naturally dual to each other.\\
\end{th}
For detail notations and results, please refer to later sections of this paper.\\\\
Our work essentially indicate that Batyrev mirror construction, which was proposed purely from toric geometry stand point, can also be understood and justified by the SYZ mirror construction. This should give us greater confidence on SYZ mirror conjecture for general Calabi-Yau manifolds.\\\\ 
This paper is a sequel of a series of papers \cite{lag1, lag2, lag3} on Lagrangian torus fibration of quintic Calabi-Yau hypersurface in ${\bf P^4}$. In \cite{lag1}, we described a very simple and natural construction via gradient flow, which will in principle be able to produce Lagrangian torus fibration for general Calabi-Yau hypersurface in toric variety. For simplicity, we described the case of Fermat type quintic Calabi-Yau threefold familly $\{X_\psi\}$ in ${\bf CP^4}$\\
\[
p_{\psi}=\sum_1^5 z_k^5 - 5\psi \prod_{k=1}^5 z_k=0
\]\\
near the large complex limit $X_{\infty}$\\
\[
p_{\infty}=\prod_{k=1}^5 z_k=0.
\]\\
in great detail. Most of the essential features of the general case already show up there. We also discussed expected special Lagrangian torus fibration structure, especially we computed the monodromy transformations of the expected special Lagrangian fibration and also discussed the expected singular fibre structure implied by monodromy information in this case. Then we compared our Lagrangian fibration constructed via gradient flow with the expected special Lagrangian fibration. Finally, we discussed its relevance to mirror construction for Calabi-Yau hypersurface in toric variety.\\\\
In order to apply our gradient flow construction effectively, it is necessary to address several technical problems concerning to the flow. First of all, the gradient flow in our situation is highly non-conventional gradient flow. The critical points of our function are usually highly degenerate and often non-isolated. Worst of all our function is not even smooth (it has infinities along some subvariety). In \cite{lag2} we discussed the local behavior of our gradient flow near critical points and infinities of our function, also the dependence on the metric to make sure that they behave the way we want.\\\\
Secondly, SYZ construction requires special Lagrangian fibration. The corresponding singular locus is expected of codimension 2. The Lagrangian torus fibration naturally produced by the gradient flow and the natural Lagrangian torus fibration of the large complex limit has singular locus of codimension 1. As described in \cite{lag1}, it is not hard to deform our Lagrangian fibration to a $C^\infty$ non-Lagrangian fibration with the same topological structure as the expected $C^\infty$ special Lagrangian fibration. To get a Lagrangian fibration of the same topological structure is more delicate. In \cite{lag2}, we are able to squeeze the codimension 1 singular locus to codimension 2 by symplectic geometry technique to get Lagrangian torus fibration with the same topological structure as the expected special Lagrangian fibration.\\\\
Thirdly, our gradient flow approach naturally produces piecewise smooth (Lipschitz) Lagrangian fibration. In our opinion, for symplectic topological aspect of SYZ construction, one should consider $C^\infty$ Lagrangian fibration. This is indicated by somewhat surprising fact, discussed in \cite{lag2}, that a general piecewise smooth (Lipschitz) Lagrangian fibration can not be smoothed to a $C^\infty$ Lagrangian fibration by small perturbation. More precisely, the singular locus of a general piecewise smooth (Lipschitz) Lagrangian fibration is usually of codimension 1 while the singular locus of the corresponding $C^\infty$ Lagrangian fibration necessarily has codimension 2. In \cite{lag2}, we constructed many local examples of piecewise smooth (Lipschitz) Lagrangian fibrations that are not $C^\infty$. In \cite{smooth}, we discussed method to deform a piecewise smooth (Lipschitz) Lagrangian fibration to a $C^\infty$ Lagrangian fibration. Although we still fall short to make our Lagrangian fibration $C^\infty$ everywhere.\\\\
With all the techniques developed in \cite{lag2}, in \cite{lag3} we are able to construct Lagrangian torus fibration with codimension 2 singular locus for general quintic Calabi-Yau hypersurface in ${\bf CP^4}$. Here the main difficulty do not come from the gradient flow. As we mentioned before that the gradient flow will always give us a Lagrangian fibration, even for more general situation. But as we learn from the Fermat type quintic case, Lagrangian fibration so constructed is usually not of the topological type of the expected $C^{\infty}$ Lagrangian (or special Lagrangian) torus fibration. We need to modify the fibration. Unlike the case of Fermat type quintic Calabi-Yau, where we can always modify the Lagrangian fibration to the expected topological type due to its nicely behaved singular locus. It is hopeless to do the same for general quintic Calabi-Yau also due to generally badly behaved singular locus. It turns out that for general quintic {\bf near} the large complex limit (in certain sense), singular locus will be much better behaved than general quintic Calabi-Yau. And we can do the same thing as in the fermat type Calabi-Yau case to get the desired Lagrangian fibration. In \cite{lag3} we also constructed Lagrangian torus fibration for the mirror of quintic Calabi-Yau.\\\\
With all these construction of Lagrangian torus fibrations of general quintic Calabi-Yau hypersurfaces and their mirror, we are able to prove the symplectic topological SYZ mirror conjecture for generic quintic Calabi-Yau hypersurface. In \cite{lag3} we also discussed a very simple computation of monodromy of Lagrangian torus fibration and singular fibres that will appear in generic Lagrangian torus fibration. These discussion in principle give us a general construction of mirror Calabi-Yau from a given generic Lagrangian torus fibration of a Calabi-Yau manifold.\\\\
In this paper, we generalize our work in \cite{lag3} to the case of general Calabi-Yau hypersurface in toric variety with respect to a reflexive polyhedron, which is exactly the situation of the Batyrev dual polyhedron mirror construction. Compared to the quintic case, in general toric case, usually both the \k moduli and the complex moduli of a Calabi-Yau hypersurface are non-trivial. The construction of the Lagrangian torus fibration has to depend on both the \k form and the complex structure of the Calabi-Yau hypersurface. We also need the most general monomial-divisor map to carry out the discussion of symplectic topological SYZ mirror construction for general Calabi-Yau hypersurface in toric variety.\\\\
In section 2, we start by reviewing basic facts from toric geometry necessary for our discussion. \k moduli of the Calabi-Yau hypersurface is discussed in section 3. Our symplectic topological SYZ construction is based on the monomial-divisor mirror map. The monomial-divisor mirror map of Aspinwall, Greene and Morrison was first constructed in \cite{AGM} with certain restriction. In sections 4 and 5, using slicing theorem, we generalize their construction of monomial-divisor mirror map to the case of Calabi-Yau hypersurface in toric variety with respect to general reflexive polyhedron. We also pin down the suitable toric compactification of the complex moduli near the large complex limit.\\\\
In section 6, we discuss the construction of the Lagrangian torus fibration of Calabi-Yau hypersurfaces in toric variety that depends on both the \k form and the complex structure of the Calabi-Yau hypersurface. At first, the gradient flow produce a Lagrangian torus fibration with codimension 1 singular locus for general Calabi-Yau hypersurface. With the miraculous fact that near the large complex and the large radius limit the singular locus behave like a fattening of a suitable graph, using techniques developed in \cite{lag2}, we are able to modify our construction to get a Lagrangian torus fibration with codimension 2 singular locus.\\\\
The symplectic topological SYZ mirror conjecture is proved in section 7 and 8. There are two major ingradients, the identification of bases (topologically $S^3$) and the duality relation of the fibres of the Lagrangian torus fibrations of the Calabi-Yau hypersurface of the toric variety and its mirror. The bases of the two fibrations are naturally boundaries of the two convex polyhedrons. Their identification can be understood purely from the combinatoric properties of the two convex polyhedrons. This identification is discussed in section 7. In section 8, we discuss the duality relation of the fibres of the Lagrangian torus fibrations of the Calabi-Yau hypersurface of the toric variety and its mirror.\\\\ 
As mentioned in the previous remark, the original SYZ conjecture is not enough to construct the mirror manifold completely. The main problem is the lack of knoledge of singular locus, singular fibres and duality of singular fibres. With our construction of Lagrangian torus fibration in \cite{lag1, lag2, lag3} and this paper, we have much better knowledge of the structure of singular locus, generic singular fibres and how singular fibres are dual to each other. This detailed knowledge enable us to precisely formulate the SYZ construction in section 9, which will enable us to construct mirror manifold topologically if a SYZ Lagrangian torus fibration of a Calabi-Yau manifold is given.\\\\
{\bf Note:} Partly based on our previous monodromy compuation in \cite{lag1}, at around the same time as our paper \cite{lag3}, in his intersting work \cite{Gross3}, M. Gross was able to use the classical knowledge of resolution of singularity of ${\bf C}^3/({\bf Z}_5\times{\bf Z}_5)$ and monodromy information to figure out the {\bf expected} fibration structure of certain quintic. This is in the same line with the idea in our previous paper \cite{lag1} to use monodromy information to figure out the {\bf expected} fibration. Gross went one significant step further to compute the cubic form etc. of the constructed topological manifold and use theorems of C.T.C. Wall to verify that the manifold consructed is indeed diffeomorphic to quintic, which is very different from our approach, therefore constructing a torus fibration for certain quintic Calabi-Yau (not necessarily Lagrangian fibration). In our work \cite{lag3}, using gradient flow approach, we were able to naturally construct Lagrangian torus fibration for general quintic Calabi-Yau and its mirror. Then we show that there is a natural duality relation between the fibres, therefore proving the symplectic SYZ completely for quintic Calabi-Yau.\\\\

\se{Background in toric geometry}
In this section, we will first review some basic constructions in toric geometry, especially the toric variety associated with a reflexive polyhedron. Then we will discuss possible singularities of simplicial toric varieties. For our gradient flow method to work effectively, this knowledge of singularities is very crucial. Finally, we will discuss some facts necessary for getting toric quotients and performing toric blow up.\\\\
\subsection{Basic}
Consider a rank $r$ lattice $N\cong {\bf Z}^r$ and its dual lattice $M=\hom_{\bf Z}(N,{\bf Z})$ with the natural pairing:\\
\[
\langle\ ,\ \rangle: M\times N \rightarrow {\bf Z}.
\]
A {\bf fan} in $N$ is a nonempty collection $\Sigma$ of strongly convex ($\sigma \cap (-\sigma)=\{0\}$) rational polyhedral cones in $N$, such that subcones of a polyhedral cone $\sigma\in \Sigma$ is also in $\Sigma$.\\\\
For a cone $\sigma\in \Sigma$, the dual cone is defined as\\
\[
\sigma^{\vee} =\{m\in M|\langle m,n\rangle \geq 0 \ {\rm for\ all}\ n\in \sigma\}.
\]
Let $U_\sigma = {\rm Spec}({\bf C}[\sigma^\vee])$. An inclusion $\sigma \hookrightarrow \sigma '$ induce an inclusion $U_\sigma \hookrightarrow U_{\sigma '}$. Gluing $U_\sigma$ together for $\sigma \in \Sigma$, we get the toric variety $P_\Sigma$. In particular, $T=U_{\{0\}} = {\rm Spec}({\bf C}[M])\cong N_{\bf C}/N=N\otimes_{\bf Z} {\bf C}^*$ is a complex $r$-torus. There is a natural action of $T=U_{\{0\}}$ on $P_\Sigma$ and $U_{\{0\}}$ is the unique Zariski dense orbit in $P_\Sigma$.\\\\
Let's fix $\sigma \in \Sigma$ and try to understand one affine piece $U_{\sigma}$ better. Now the problem can be rephrased as given an $r$ dimensional rational convex polyhedral cone $\sigma^{\vee}\subset M$, discuss the structure of $U_\sigma = {\rm Spec}({\bf C}[\sigma^\vee])$. For any cone $\tau\subset M$, denote $X_{\tau}={\rm Spec}({\bf C}[\tau])$. (Relate to the old notation $X_{\sigma^{\vee}} = U_{\sigma}$.)\\\\
For a subcone $\tau\subset \sigma^{\vee}$, There is a seemingly unnatural semi-group map $r: \sigma^{\vee}\rightarrow \tau$ ($r|_\tau = id$, $r|_{\sigma^\vee\backslash \tau}=0$) that induce a very natural embedding $X_{\tau}\hookrightarrow X_{\sigma^{\vee}} =U_{\sigma}$. $r$ is just the corresponding restriction map. On the other hand, the natural embedding $\tau\subset \sigma^{\vee}$ induces a seemingly rather unnatural projection $X_{\sigma^{\vee}} =U_{\sigma}\rightarrow X_{\tau}$. (This projection can be understood as every torus in $X_{\sigma^{\vee}} =U_{\sigma}$ is projected to the largest adjacent torus in $X_{\tau}$. In particular, fibre dimension could jump.)  Let $\overline{\tau}$ denote the span of $\tau$, and $T_{\tau}= X_{\overline{\tau}}$. We have\\
\begin{center}
\setlength{\unitlength}{1pt}
\begin{picture}(100,50)(0,0)
\put(-50,0){
\put(30,40){\makebox(0,0){$X_{\tau}$}}
\put(80,40){\makebox(0,0){$X_{\sigma^{\vee}}$}}
\put(30,10){\makebox(0,0){$T_{\tau}$}}
\put(80,10){\makebox(0,0){$T_{\sigma^{\vee}}$}}
\put(43,42){\oval(6,4)[l]}
\put(43,40){\vector(1,0){23}}
\put(30,20){\vector(0,1){12}}
\put(28,20){\oval(4,6)[b]}
\put(80,20){\vector(0,1){12}}
\put(78,20){\oval(4,6)[b]}
\put(38,17){\vector(2,1){30}}}

\put(50,0){
\put(30,40){\makebox(0,0){$X_{\tau}$}}
\put(80,40){\makebox(0,0){$X_{\sigma^{\vee}}$}}
\put(30,10){\makebox(0,0){$T_{\tau}$}}
\put(80,10){\makebox(0,0){$T_{\sigma^{\vee}}$}}
\put(66,40){\vector(-1,0){26}}
\put(30,20){\vector(0,1){12}}
\put(28,20){\oval(4,6)[b]}
\put(66,10){\vector(-1,0){26}}
\put(80,20){\vector(0,1){12}}
\put(78,20){\oval(4,6)[b]}}
\end{picture}
\end{center}
Let $\Sigma(k)$ denote the collection of $k$-dimensional cones in the fan $\Sigma$, then for any $\sigma\in\Sigma(1)$, $X_{\sigma^\perp}\hookrightarrow U_{\sigma}\hookrightarrow P_\Sigma $ give us a $T$-invariant irreducible Weil divisor $D_{\sigma}$. Let $e_\sigma$ denote the unique primitive element in $\sigma$. Consider $M$ as the set of $T$-invariant meromorphic functions on $P_\Sigma$, then for $m\in M$, $\langle m,e_\sigma\rangle$ is the vanishing order of $m$ along $D_\sigma$. This important observation motivates the representation of line bundle by piecewise linear function on the fan. The most important line bundle is the anti-canonical bundle ${\cal O}(D_\Sigma)$, where\\
\[
D_\Sigma = \bigcup_{\sigma\in \Sigma(1)} D_\sigma
\]\\
${\cal O}(D_\Sigma)$ can be characterized by the piecewise linear function $p_{\Sigma}$ that satisfy $p_\Sigma(e_\sigma)=-1$ for any $\sigma\in \Sigma(1)$.\\\\
Given a fan $\Sigma$, there are two naturally associated polyhedrons. The Newton polyhedron of ${\cal O}(D_\Sigma)$ defined as\\
\[
\Delta_\Sigma = \{m\in M|\langle m,n\rangle \geq p_{\Sigma}(n),\ \ {\rm for\ any\ } n\in N\}
\]\\
The second polyhedron $\Delta_\Sigma^\vee$ is the convex hull of $e_\sigma$ for $\sigma\in \Sigma(1)$. Since $p_{\Sigma}(e_\sigma)=-1$ for $\sigma\in \Sigma(1)$, if we think of $M$, $N$ as real vector spaces, the two real polyhedrons are naturally dual to each other. Namely\\
\[
\Delta_\Sigma^\vee = \{n\in N|\langle m,n\rangle \geq -1,\ \ {\rm for\ any\ } m\in \Delta_\Sigma \}
\]\\
In our situation, we think of $M$, $N$ as two lattices. $\Delta_\Sigma$ and $\Delta_\Sigma^\vee$ are two integral polyhedrons. In general, we only have $\Delta_\Sigma$ is dual of $\Delta_\Sigma^\vee$, not vise versa. When the two integral polyhedrons are dual to each other, we call $\Delta_\Sigma$ {\bf reflexive}.\\\\
A dimemsion $r$ cone $\sigma$ is called {\bf simplicial} cone if it has exactly $r+1$ 1-dimensional subcones, if further more the primitive elements in those 1-dimensional subcones generate $\sigma$, then $\sigma$ is called {\bf primitive} simplicial cone.\\\\
We can also introduce volume $v_\sigma$ of a simplicial cone $\sigma$ to be the volume of the parallelgram generated by generators of 1-dimensional subcones of $\sigma$. Then a simplicial cone $\sigma$ is primative if and only if $v_\sigma=1$.\\\\
$P_\Sigma$ is usually not smooth, since $\sigma \in \Sigma$ may not be simplicial. To make $P_\Sigma$ smoother, we need to subdivide the fan $\Sigma$. We usually do not want the polyhedron $\Delta_\Sigma$ to be affected by subdivision. The best way to archieve this is to consider subdivision $\Sigma'$ with $\Sigma'(1)$ consists of only the boundary integral points of $\Delta_\Sigma^\vee$. This kind of subdivision is called {\bf crepant} subdivision.\\\\
Clearly, through crepant subdivision, we can always make the cones in $\Sigma$ simplicial. It is usually impossible to make them primitive in general only by crepant subdivision. Recall that $P_\Sigma$ is smooth if and only if each $\sigma\in \Sigma$ is primitive simplicial cone.\\\\
\subsection{Toric varieties associated with a reflexive polyhedron}
A good way to avoid many unpleasant features is to start with a reflexive polyhedron $\Delta\in M$ (when both $\Delta$ and $\Delta^\vee$ are integral). Fan $\Sigma$ should be constructed in two steps. First construct all the normal cones, for a face $\alpha$ of the polyhedron $\Delta$, define the normal cone of $\alpha$\\
\[
\sigma_\alpha :=\{n\in N|\langle m',n\rangle \leq \langle m,n\rangle \ \ {\rm for\ all\ } m'\in \alpha,\ m\in \Delta \}
\]\\
When taking $\alpha$ to be vertices of $\Delta$, we get the top dimension cones. The collection of all the normal cones $\sigma_\alpha$ forms the anti-canonical fan $\Sigma_\Delta$ of $\Delta$. In the second step, through crepant subdivision, we subdivide these cones to get a simplicial fan $\Sigma$ with as many cones as possible. In this way, $\Sigma(1)$ still consists of only the boundary integral points of $\Delta^\vee$. We will denote divisors come from the first step principle divisors and denote by $\Sigma(1)_p$, divisors come from the second step exceptional divisors and denote by $\Sigma(1)_e$.\\\\ 
Notice the second step is not unique. In this way, we are getting a collection of toric varieties $P_\Sigma$ with common $\Sigma(1)$ that is associated with reflexive polyhedron $\Delta$. They differ only by codimension two birational modification.\\\\
A fan $\Sigma$ so constructed is a maximal simplicial fan that is compatible with $\Delta$, where $\Sigma$ is called {\bf compatible} with $\Delta$ if $\Sigma$ is a crepant subdivision of the anti-canonical fan $\Sigma_\Delta$. In this case, $P_\Sigma$ is a crepant resolution of $P_{\Sigma_\Delta}$.\\\\
In our point of view, for the fan $\Sigma$, the set of 1-dimensional cones $\Sigma(1)$ is more fundamental. Let $e_\sigma$ be the primitive generator of $\sigma\in \Sigma(1)$, then $\sigma \rightarrow e_\sigma$ defines an identification of $\Sigma(1)$ with integral points in $\partial \Delta^\vee$. Namely, $\Sigma(1)$ contains the same information as $\Delta$.\\\\
\subsection{Subtoric varieties}
Recall for $\sigma\in \Sigma$, the restriction map $r: \sigma^\vee \rightarrow \sigma^\perp$ induces the embedding $T_\sigma= {\rm Spec}({\bf C}[\sigma^\perp])\hookrightarrow U_\sigma= {\rm Spec}({\bf C}[\sigma^\vee])\hookrightarrow P_\Sigma$. $\dim_{\bf C} T_\sigma + \dim \sigma = r$. (Notice here we changed notation of $T_\sigma$ a little bit.) In this way, we have a decomposition of $P_\Sigma$ into disjoint union of complex tori.\\
\[
P_\Sigma = \bigcup_{\sigma\in \Sigma} T_\sigma
\]\\
$T_\tau$ is adjacent to $T_\sigma$ (namely, $T_\tau\subset \overline{T_\sigma}$) if and only if $\sigma$ is a subcone of $\tau$. For $\sigma \in \Sigma$, consider projection $\pi_\sigma: N\rightarrow N_\sigma:=N/{\overline{\sigma}}$. Then $N_\sigma$ is naturally dual to $M_\sigma:=\sigma^\perp$. We can define a new fan\\
\[
\Sigma_\sigma = \{\pi_\sigma(\tau)\subset N_\sigma|\sigma\subset \tau \in \Sigma\}
\]\\
The corresponding toric variety $P_{\Sigma_\sigma}$ is the closure of $T_\sigma$ in $P_\Sigma$.\\\\
On the other hand, the natural injection $\sigma^\perp \rightarrow \sigma^\vee$ induces the fibration $U_\sigma= {\rm Spec}({\bf C}[\sigma^\vee]) \rightarrow T_\sigma= {\rm Spec}({\bf C}[\sigma^\perp])$.\\
\begin{prop}
\label{be}
The fibres of the fibration $U_\sigma= {\rm Spec}({\bf C}[\sigma^\vee]) \rightarrow T_\sigma= {\rm Spec}({\bf C}[\sigma^\perp])$ are naturally affine toric variety with fan $\Sigma^\sigma$ consists of subcones of $\sigma \subset \bar{\sigma}={\rm Span}(\sigma)$.\\
\end{prop}
It is interesting to understand the restriction of a line bundle on $P_\Sigma$ to a subvariety $P_{\Sigma_\sigma}$. As we know, a piecewise linear function $p$ on $N$ with respect to the fan $\Sigma$ will give us a meromorphic section $s_p$ of the corresponding line bundle $L_p$. $s_p$ is non-vanishing (vanishing, has a pole) when restrict to the torus $T_\sigma$ if and only if $p$ restrict to the interior of $\sigma$ is zero (negative, positive). In particular, $s_p$ is holomorphic if and only if $p\leq 0$. To make restriction to $P_{\Sigma_\sigma}$ meaningful, we would want section $s_p$ of $L_p$ to be holomorphic and non-vanishing along $T_\sigma$. Namely, we want $p$ restrict to zero on $\sigma$. Since $p$ is convex, it is easy to modify $p$ by a linear function $m\in M$ to ensure $p\leq 0$ and $p|_\sigma=0$. Under these conditions, $p$ will naturally induce a piecewise linear convex function $p_\sigma$ on $N_\sigma$ with respect to fan $\Sigma_\sigma$. Now assume cones in $\Sigma$ are all simplicial cones.\\\\
\subsection{Singularities of a simplicial toric variety}
As we know, $P_\Sigma$ is an union of complex torus (orbits of $T$ action) that is parametrized by cones in the fan $\Sigma$. In particular, $r$ dimensional cones $\sigma\in \Sigma(r)$ give rise to zero dimensional orbits (fixed points of the $T$ action) $q_\sigma$. On the other hand, an $r$ dimensional cone $\sigma\in \Sigma(r)$ also gives rise to an affine toric subvariety $U_\sigma\in P_\Sigma$, which is naturally a neighborhood of $q_\sigma$. All this neighborhoods $\{U_\sigma\}_{\sigma\in \Sigma(r)}$ form an open covering of $P_\Sigma$.\\\\
\begin{prop}
\label{bf}
Assume $e_1,e_2,e_3$ are the generators of a 3-dimensional simplicial cone $\sigma$ in a lattice $M$ and $\{0,e_1,e_2,e_3\}= \sigma\cap\Delta$ for reflexive polyhedron $\Delta\subset M$. Then $\sigma$ is primative.\\
\end{prop}
{\bf Proof:} Assume that $\sigma$ is not primative, there exist a prime number $p$ and non-negative integers $a,b,c$ such that\\
\[
m= \frac{a}{p}e_1 + \frac{b}{p}e_2 + \frac{c}{p}e_3\not=0
\]\\
is in $\sigma$. Without loose of generality, we may assume $c=1$ and $0\leq a,b\leq p-1$. Then $1\leq a+b+c\leq 2p-1$.\\\\
Since $\{0,e_1,e_2,e_3\}=\sigma\cap\Delta$, by theorem \ref{bb}, $e_1,e_2,e_3$ are all in one codimension 1 face of reflexive polyhedron $\Delta\subset M$, there exist an element $n$ in the dual lattice of $M$ such that $n|_\Delta\leq 1$ and $\langle e_k,n\rangle=1$ for $k=1,2,3$. Since $m\in M$,\\
\[
\langle m,n\rangle=\frac{a+b+c}{p}
\]\\
should be an integer. Therefore\\
\[
a+b+c=p
\]\\
$m$ as a convex combination of $e_1,e_2,e_3$ is in $\sigma\cap\Delta$. A contradiction!
\begin{flushright} $\Box$ \end{flushright}
{\bf Remark:} This theorem is equivalent to the fact that in any rank two lattice, a simplex with integral vertex that do not contain integral points in inside and edges is always the standard simplex with area $\frac{1}{2}$.\\\\
By proposition \ref{be} and \ref{bf}, we have\\
\begin{co}
\label{ba}
If $\Sigma$ is a maximal simplicial fan on a 4-dimensional lattice $N$ that is compatible with a reflexive polyhedron $\Delta\subset M$, Then $P_\Sigma$ is smooth away from the zero dimension subtorus (fixed points under the action of $T$) of $P_\Sigma$.
\end{co}
\begin{flushright} $\Box$ \end{flushright}
\begin{co}
\label{bb}
If $\Sigma$ is a maximal simplicial fan on a 4-dimensional lattice $N$ that is compatible with a reflexive polyhedron $\Delta\subset M$, then a generic section of the anti-canonical bundle $L_\Delta$ on $P_\Sigma$ cut out a smooth 3-dimensional Calabi-Yau manifold.
\end{co}
\begin{flushright} $\Box$ \end{flushright}
Or more generally\\
\begin{co}
If $\Sigma$ is a maximal simplicial fan on a lattice $N$ that is compatible with a reflexive polyhedron $\Delta\subset M$, Then $P_\Sigma$ is smooth away from the codimension 4 subtorus of $P_\Sigma$. A generic section of the anti-canonical bundle $L_\Delta$ on $P_\Sigma$ cut out a Calabi-Yau hypersurface with codimension 4 singularities.
\end{co}
\begin{flushright} $\Box$ \end{flushright}
\subsection{Toric action of a torus on a toric variety}
As we know, $M$ can be viewed as monomial functions on $T= {\rm Spec}({\bf C}[M])\cong N_{\bf C}/N=N\otimes_{\bf Z} {\bf C}^*$. Given a subspace $W\subset N$, $W^\perp\subset M$ can be viewed as $T_W=W_{\bf C}/W=W\otimes_{\bf Z} {\bf C}^*$ invariant functions on $T$, or on another word, functions on $T/T_W$.\\\\
For $\sigma\in \Sigma$, the action of $T_W$ on $T$ can be naturally extended to an action on $T_\sigma$, that is indicated by $W_\sigma\subset N_\sigma = N/\bar{\sigma}$. The stablizer of the action is corresponding to $W\cap \bar{\sigma}$. We have the exact sequence\\
\[
0\rightarrow T_{W\cap \bar{\sigma}}\rightarrow T_W \rightarrow T_{W_\sigma} \rightarrow 0
\]\\
$W^\perp\cap \sigma^\perp \subset \sigma^\perp$ can be viewed as $T_{W_\sigma}=W_\sigma\otimes_{\bf Z} {\bf C}^*$ invariant functions on $T_\sigma$, or on another word, functions on $T_\sigma/T_{W_\sigma}$.\\
\begin{lm}
\label{bc}
$T_W$ acts freely on $T_\sigma$ if and only if $\bar{\sigma}\cap W= \{0\}$.
\end{lm}
\begin{flushright} $\Box$ \end{flushright}
\subsection{Blow up}
Given two rational strongly convex polyhedron cones $\tau\subset\sigma$ in $N$ ($\tau$ is not necessaryly a {\bf subcone} of $\sigma$), then\\
\begin{prop}
$\tau^\perp \cap \sigma^\vee$ is a subface of $\sigma^\vee$, and the following diagram commutes\\
\begin{center}
\setlength{\unitlength}{1pt}
\begin{picture}(100,75)(0,0)
\put(-50,0){
\put(80,70){\makebox(0,0){$\sigma_\tau^\perp$}}
\put(30,40){\makebox(0,0){$\tau^\perp$}}
\put(85,40){\makebox(0,0){$\tau^\perp \cap \sigma^{\vee}$}}
\put(30,10){\makebox(0,0){$\tau^\vee$}}
\put(80,10){\makebox(0,0){$\sigma^{\vee}$}}
\put(60,40){\oval(6,4)[r]}
\put(60,38){\vector(-1,0){20}}
\put(63,10){\oval(6,4)[r]}
\put(63,8){\vector(-1,0){23}}
\put(30,17){\vector(0,1){15}}
\put(80,17){\vector(0,1){15}}
\put(80,47){\vector(0,1){15}}
\put(68,62){\vector(-2,-1){30}}}

\put(50,0){
\put(80,70){\makebox(0,0){$T_{\sigma_\tau}$}}
\put(30,40){\makebox(0,0){$T_{\tau}$}}
\put(85,40){\makebox(0,0){$X_{\tau^\perp \cap \sigma^{\vee}}$}}
\put(30,10){\makebox(0,0){$U_{\tau}$}}
\put(80,10){\makebox(0,0){$U_{\sigma}$}}
\put(38,47){\vector(2,1){30}}
\put(80,59){\vector(0,-1){12}}
\put(78,59){\oval(4,6)[t]}
\put(40,40){\vector(1,0){23}}
\put(30,29){\vector(0,-1){12}}
\put(28,29){\oval(4,6)[t]}
\put(40,10){\vector(1,0){26}}
\put(80,29){\vector(0,-1){12}}
\put(78,29){\oval(4,6)[t]}}
\end{picture}
\end{center}
where $\sigma_\tau$ is the smallest subcone of $\sigma$ that contains $\tau$ and clearly, $\sigma_\tau^\perp = {\rm Span}( \tau^\perp \cap \sigma^{\vee})$.
\end{prop}
\begin{flushright} $\Box$ \end{flushright}
\begin{co}
\label{bd}
If $\dim \tau =\dim \sigma$, then $\sigma_\tau =\sigma$, $T_\tau\rightarrow T_\sigma$ is an isomorphism.
\end{co}
\begin{flushright} $\Box$ \end{flushright}
This proposition implies that $U_{\tau}$ is blow up from $U_{\sigma}$ and also indicates that how subtorus in $U_{\tau}$ is blow down to subtorus of $U_{\sigma}$.\\\\

\se{\k moduli of Calabi-Yau hypersurface in toric variety}
It is interesting to consider Calabi-Yau hypersurface in a toric variety $P_\Sigma$ that comes from a reflexive polyhedron $\Delta\subset M$. As we know, anti-canonical bundle $L_{\Delta}$ is semi-ample. A generic section $s$ of $L_{\Delta}$ cut out a Calabi-Yau hypersurface $X=X_s\subset P_\Sigma$. Integral points in $\Delta$ naturally give us $T$-invariant sections of $L_{\Delta}$, $s$ can be taken to be generic linear combination of these $T$-invariant sections. For the purpose of mirror symmetry, it is important to understand the complex moduli and \k cone of these Calabi-Yau hypersurfaces. Before describing the \k moduli of Calabi-Yau hypersurface, first we need to discuss the \k cone of a toric variety.\\\\
\subsection{\k cone of a toric variety}
It is interesting to understand the \k cone of $P_\Sigma$. Namely the cone of convex piecewise linear functions modulo linear functions. Recall that $m\in M$ (linear function on $N$) corresponds to $T$ invariant meromorphic functions. An integral piecewise linear function $p$ corresponds to a meromorphic section (up to multiplication by a non-zero constant) of the corresponding line bundle, those $p\leq 0$ correspond to holomorphic sections. In particular, $m=0\leq 0$ corresponds to constant functions that are the only holomorphic functions on $P_\Sigma$. More precisely, an integral piecewise linear function corresponds to a divisor. (A divisor is equivalent to a line bundle together with a non-trivial meromorphic section up to non-zero constant multiple.)\\\\
A piecewise linear function $p$ is called {\bf compatible} with a fan $\Sigma$ if the restriction of $p$ to each cone $\sigma\in \Sigma$ is linear. Fix value of $p$ on $\Sigma(1)$, there are many ways to extend $p$ to piecewise linear function on $N$. For any simplicial fan $\Sigma$ that is compatible with the Newton polyhedron $\Delta_\Sigma$ (determined by $\Sigma(1)$), there is an unique extension of $p$ as a piecewise linear function on $N$ that is compatible with $\Sigma$.\\\\
Conversely, a piecewise linear function $p$ will canonically determine a fan $\Sigma_p$ such that $\Sigma_p$ is compatible with $p$ and any other $\Sigma$ compatible with $p$ is a subdivision of $\Sigma_p$. A (convex) piecewise linear function $p$ is called strongly (convex) piecewise linear with respect to a fan $\Sigma$ compatible with $p$ if $\Sigma = \Sigma_p$, we will also say $p$ is strongly (convex) piecewise linear without mentioning the fan $\Sigma$ if $\Sigma=\Sigma_p$ is simplicial.\\\\
Let $\Sigma$ be a simplicial cone compatible with the Newton polyhedron $\Delta$, then the classical \k cone is just the set of (strongly) convex piecewise liear function $p$ (modulo linear functions $M$) that is compatible with $\Sigma$. For our purpose, we are more interested in the movable cone, which is solely determined by the polyhedron $\Delta$. For a piecewise linear function $p$, let $p^1=p|_{\Sigma(1)}$. We call $p^1$ the 1-skeleton of $p$. To define the movable cone, in general we also use $p^1$ to denote a general piecewise linear finction on $\Sigma(1)$. Then the movable cone\\
\[
K_\Delta = \{p^1|{\rm for\ any}\ \sigma\in \Sigma(1),\ {\rm exist}\ m_\sigma\in M,\ {\rm s.t.}\ m_\sigma\geq p^1, \ m_\sigma|_\sigma=p^1|_\sigma\}/M
\]\\
It is easy to see that $K_\Delta$ is a convex cone. For any $p^1\in K_\Delta$, we can define\\
\[
p(n)=\min\{\langle m,n\rangle | m\geq p^1,\ m\in M\}
\]\\
Clearly, $p$ is a convex piecewise linear extension of $p^1$, and $p^1$ is the 1-skeleton of $p$. Namely, if the 1-skeleton $p^1$ of $p$ is in $K_{\Delta}$, then $p$ is completely determined by $p^1$ and $p$ is convex. For generic $p^1\in K_{\Delta}$, $p$ is strongly convex, $\Sigma_p$ is a simplicial fan. According to type of $\Sigma_p$, the movable cone $K_\Delta$ is divided into subcones that are the \k cones of different birational equivalent models.\\\\
For $p^1\leq 0$, a piecewise linear function on $\Sigma(1)$, for each $\sigma\in \Sigma(1)$, we can define halfspace\\
\[
H_\sigma =\{m\in M| \langle m,e_\sigma\rangle \geq p^1(e_\sigma)\}
\]\\
\[
\Delta_p =\bigcap_{\sigma\in \Sigma(1)} H_\sigma
\]\\
is a convex polyhedron in $M$ containing origion. $p^1\in K_\Delta$ if and only if for each $\sigma\in \Sigma(1)$, $\Delta_p$ touch boundary of $H_\sigma$.\\\\
In particular, when $\Delta$ is reflexive the anti-canonical divisor corresponds to $p_\Delta^1(e_\sigma)=-1$ for $\sigma\in \Sigma(1)$. The corresponding convex polyhedron $\Delta_{p_\Delta}$ is exactly $\Delta$. $H_\sigma$ will touch $\Delta$ in dimension $r-1$ when $\sigma\in\Sigma(1)_p$, and will touch in lower dimension when $\sigma\in\Sigma(1)_e$. Therefore anti-canonical class is in $K_\Delta$. Further more, $p_\Delta$ is compatible with any simplicial fan $\Sigma$ that is compatible with $\Delta$. Therefore, anti-canonical class is in the \k cone of every birational equivalent model that dominate the anti-canonical model. Namely in the intersection of all the \k cones of the birational equivalent models that dominate the anti-canonical model.\\\\
{\bf Remark:} In general, for $p^1\in K_\Delta$, apriorily, the corresponding $p$ is not necessary compatible with the Newton polyhedron $\Delta$. This corresponds to the question, weither the corresponding toric variety will dominate anti-canonical model or in another word, weither the anti-canonical bundle is semi-ample. It turns out the anti-canonical bundle is always semi-ample.\\\\
\begin{th}
\label{ca}
For strongly convex $p^1\in K_\Delta$, the corresponding $p$ is always compatible with the Newton polyhedron $\Delta$. Equivalently, $P_{\Sigma_p}$ dominate the anti-canonical model.\\
\end{th}
This result is easyly implied by the following\\
\begin{th}
\label{cb}
For reflexive polyhedron $\Delta\in M$, if two integral points $m_1,m_2\in \Delta$ is not in the same face of $\Delta$, then $m_1 + m_2\in \Delta$.\\
\end{th}
{\bf Proof:} Recall that corresponging to the reflexive polyhedron $\Delta$ there is an integer valued non-positive convex function $p_\Delta$, such that\\
\[
\Delta = \{m\in M|p_\Delta(m)\geq -1\}.
\]
\[
\partial \Delta = \{m\in M|p_\Delta(m) = -1\}.
\]
Hence $p_\Delta(m_1) = p_\Delta(m_2) =-1$. Convexity of $p_\Delta$ implies that\\
\[
p_\Delta(m_1 + m_2) \geq p_\Delta(m_1) + p_\Delta(m_2)=-2.
\]
Equality hold if and only if $m_1,m_2$ are in the same face of $\Delta$. Therefore, when $m_1,m_2$ are not in the same face of $\Delta$, we have $p_\Delta(m_1 + m_2)=-1$ and $m_1 + m_2$ is in $\Delta$.
\begin{flushright} $\Box$ \end{flushright}
{\bf Remark:} The theorem \ref{cb} proved more than what is needed for the theorem \ref{ca}, and is interesting in its own.\\
\begin{co}
The movable cone $K_\Delta$ of $P_\Sigma$ is divided into union of \k cones of different birational models. All of these \k cones have a common boundary point, the anti-canonical class of $P_\Sigma$.\\
\end{co}  
Given $\Sigma(1)$ consists of integral points on the boundary of $\Delta^\vee$, take a convex body $B$ containing origion (e.g. standard ball), then for any $\sigma\in \Sigma(1)$, $B$ intersects $\sigma$ at a vector $v_\sigma$. $p^1(v_\sigma)=-1$ will give us $p^1\in K_\Delta$.\\\\
\subsection{Restriction to a Calabi-Yau hypersurface}
The toric part of the \k cone of $X$ is the restriction of the \k cone of $P_\Sigma$ to $X\subset P_\Sigma$. Just like in the toric case, it is more natural to consider movable cone of $X$. Recall that the \k cone of $P_\Sigma$ corresponds to convex piecewise linear function $p$ on $N$ with respect to $\Sigma$. Function $p$ are determined by its restriction $p^1$ of $p$ to $\sigma \in \Sigma(1)$, whose primitive elements $e_\sigma$ are vertices of $\Delta^\vee$. The union of all the \k cones of $P_\Sigma$ for $\Sigma$ compatible with $\Delta$ forms the movable cone. Anti-canonical class is determined by $p_\Delta$ that satisfies $p_\Delta(e_\sigma) = -1$ for $\sigma\in \Sigma(1)$. $p_\Delta$ also determine a section of $L_\Delta$ that corresponds to $0\in \Delta$, which vanishes to order one on each divisor $D_\sigma=P_{\Sigma_\sigma} \subset P_\Sigma$ for $\sigma\in \Sigma(1)$. Recall to get the restriction of $L_\Delta$ to $D_\sigma$, we need to adjust $p_\Delta$ by linear functions in $M$ to get a $p_\sigma$ such that $p_\sigma|_\sigma=0$ and $p_\sigma\leq 0$. $p_\sigma$ will naturally give us a piecewise linear function (still denote by $p_\sigma$) on $N_\sigma$ with respect to the fan $\Sigma_\sigma$. It is easy to see that\\
\begin{lm}
$p_\sigma=0$ on $N_\sigma$ if and only if $e_\sigma$ is an inner point of a $r-1$ dimensional face of $\Delta^\vee$.\\
\end{lm}
\begin{co}
$L_\Delta$ restricts to ${\cal O}_{D_\sigma}$ on $D_\sigma$ if and only if $e_\sigma$ is an inner point of a $r-1$ dimensional face of $\Delta^\vee$.\\
\end{co}
\begin{co}
For a generic section $s$ of $L_\Delta$, $X_s$ does not intersect $D_\sigma$ if and only if $e_\sigma$ is an inner point of a $r-1$ dimensional face of $\Delta^\vee$.\\
\end{co}
Let $\Sigma(1)_0$ be the subset of $\Sigma(1)$ consists of $\sigma\in \Sigma(1)$ such that $e_\sigma$ is not an inner point of a $r-1$ dimensional face of $\Delta^\vee$. Let $p^0=p|_{\Sigma(1)_0}$, then the movable cone $K_\Delta^0$ of $X$ can be defined as\\
\[
K_\Delta^0 = \{p^0|{\rm for\ any}\ \sigma\in \Sigma(1)_0,\ {\rm exist}\ m_\sigma\in M,\ {\rm s.t.}\ m_\sigma\geq p^0, \ m_\sigma|_\sigma=p^0|_\sigma\}/M
\]\\
It is easy to see that $K_\Delta^0$ is a convex cone. For any $p^0\in K_\Delta^0$, we can define\\
\[
p(n)=\min\{\langle m,n\rangle | m\geq p^0,\ m\in M\}
\]\\
Clearly, $p$ is a convex piecewise linear extension of $p^0$. For generic $p^0\in K_{\Delta}^0$, $p$ is strongly convex, $\Sigma_p$ is a simplicial fan. According to type of $\Sigma_p$, the movable cone $K_\Delta^0$ is divided into subcones that are the \k cones of different birational equivalent models of $X$.\\\\
The two movable cones $K_\Delta$ and $K_\Delta^0$ are closely related. There are natural projection map $\pi: K_\Delta \rightarrow K_\Delta^0$ and inclusion map $i: K_\Delta^0 \hookrightarrow K_\Delta$. $\pi$ is defined as $p^0=p^1|_{\Sigma(1)_0}$, and $i$ is defined as the restriction of the extension $p$ of $p^0$ to $\Sigma(1)$. Clearly, $\pi\circ i = id$.\\\\

\se{Automorphism group of toric variety and slicing theorem}
The toric part of the complex moduli of $X=X_s$ is parametrized by section $s$ of the anti-canonical bundle $L_\Delta$ modulo the action of the automorphism group of $P_\Sigma$. The space of sections modulo constant multiple is a projective space. To understand the complex moduli, the key point is to understand the automorphism group of $P_\Sigma$. In this section, we will first discuss the automorphism group of $P_\Sigma$, then we will prove a slicing theorem that enable us to reduce the complex moduli near the large complex limit from a quotient by ${\rm Aut}(P_\Sigma)$ to a toric quotient by the maximal torus of ${\rm Aut}(P_\Sigma)$.\\\\ 
\subsection{Automorphism group of a toric variety} 
The automorphism group ${\rm Aut}(P_\Sigma)$ of $P_\Sigma$ is a complex Lie group. Infinitesimally, automorphism group can be decomposed into two parts, the part that keep the subtorus of $P_\Sigma$ invariant (the toric part) and the part that move those subtorus (the non-toric part). Elements in the first part are determined by the automorphisms of the big torus $T\cong U_0$. Elements in the second part (more precisely, modulo the first part) are determined by their action on each divisor $D_\sigma$ for $\sigma\in \Sigma(1)$, which can be characterized as sections of the normal bundles ${\cal O}(D_\sigma)|_{D_\sigma}$ of $D_\sigma \cong P_{\Sigma_\sigma}$ for $\sigma\in \Sigma(1)$. On the other hand, the action of non-toric part of ${\rm Aut}(P_\Sigma)$ on divisors $D_\sigma$ for $\sigma\in \Sigma(1)$ are independent. Namely, any configuration of sections of ${\cal O}(D_\sigma)|_{D_\sigma}$ for $\sigma\in \Sigma(1)$ can be realized from an infinitesimal automorphism of $P_\Sigma$. More precisely, this discription of the infinitesimal automorphism group of $P_\Sigma$ (the Lie algebra) can be expressed by the following exact sequence:\\
\begin{lm}
\[
0\rightarrow \Theta_{P_\Sigma}(-\log D) \rightarrow \Theta_{P_\Sigma} \rightarrow \bigoplus_{\sigma\in \Sigma(1)} {\cal O}(D_\sigma)|_{D_\sigma} \rightarrow 0
\]
is exact, where $D=\sum_\sigma D_\sigma$, $\Theta_{P_\Sigma}$ denote the tangent sheaf of $P_\Sigma$.\\
\end{lm}
\begin{co}
\label{db}
\[
0\rightarrow H^0(\Theta_{P_\Sigma}(-\log D)) \rightarrow H^0(\Theta_{P_\Sigma}) \rightarrow \bigoplus_{\sigma\in \Sigma(1)} H^0({\cal O}(D_\sigma)|_{D_\sigma}) \rightarrow 0
\]
is exact.\\
\end{co}
\begin{prop}
\label{dc}
\[
\Theta_{P_\Sigma}(-\log D) \cong {\cal O}_{P_\Sigma}\otimes_{\bf Z} N
\]
\[
H^0(\Theta_{P_\Sigma}(-\log D)) \cong {\bf C}\otimes_{\bf Z} N
\]\\
\end{prop}
Let $H_\sigma = \{m\in M| \langle m, e_\sigma\rangle =-1\}$, $\Delta \cap H_\sigma$ is a subface of $\Delta$ of dimension $\leq r-1$. When $\Delta \cap H_\sigma$ is of dimension $r-1$, let $\Delta_\sigma$ denotes the polyhedron consisting of interior points of $\Delta \cap H_\sigma$. When $\Delta \cap H_\sigma$ is of dimension $<r-1$, let $\Delta_\sigma = \emptyset$. $\Delta_\sigma$ can also be interpreted as the set of supporting hyperplanes of $\Delta^\vee$ that touch $\Delta^\vee$ only at $e_\sigma$.\\
\begin{prop}
For $\sigma\in \Sigma(1)$, $\Delta_\sigma$ is exactly the Newton polyhedron of the line bundle ${\cal O}(D_\sigma)|_{D_\sigma}$. Namely,\\
\[
H^0({\cal O}(D_\sigma)|_{D_\sigma}) \cong {\rm Span}_{\bf C}(\Delta_\sigma).
\]\\
\end{prop}
As we know, vector fields on $T=U_0={\rm Spec}({\bf C}[M])$ can be understood as differentials from ${\bf C}[M]$ to itself. First, let's consider differentials that respect the toric structure. Namely, those with $m\in M$ as eigenvectors. This kind of ${\bf C}$-linear maps on ${\bf C}[M]$ is differential if and only if the eigenvalue $\lambda_m$ of $m\in M$ as a function on $M$ is linear with repect to the linear structure on $M$. Hence, for each $n\in N$, we can define a differential (vector field) $v_n$ that respect the toric structure by  $v_n(m)=\langle m,n\rangle m$. This naturally gives us the embedding\\
\[
{\bf C}\otimes_{\bf Z} N \hookrightarrow H^0(\Theta_{U_0}) 
\]\\
It is easy to see that these vector fields can naturally be extended to holomorphic vector fields on $P_\Sigma$ and form all the sections of $\Theta_{P_\Sigma}(-\log D)$. These vector fields vanish to order at most 1 at each divisor $D_\sigma$ for $\sigma\in \Sigma(1)$. $v_n$ vanishes at $D_\sigma$ if and only if $n$ is a multiple of $e_\sigma$. General holomorphic vector fields on $U_0$ can be characterized as\\
\begin{prop}
\[
H^0(\Theta_{U_0}) \cong {\bf C}[M]\otimes_{\bf Z} N
\]\\
\end{prop}
In general, holomorphic vector fields on $P_\Sigma$ can be characterized by image of the restriction map $H^0(\Theta_{P_\Sigma}) \hookrightarrow H^0(\Theta_{U_0})$.\\
\begin{prop}
$mv_n$ can be extended holomorphically to $P_\Sigma$ if and only if $n$ is a multiple of $e_\sigma$ for a $\sigma\in \Sigma(1)$ and $m$ defines a supporting function of $\Delta^\vee$, whose supporting hyperplane touches $\Delta^\vee$ only at $e_\sigma$.\\
\end{prop}
Now we have a natural splitting of the exact sequence in corollary \ref{db}\\
\begin{co}
\label{de}
\[
H^0(\Theta_{P_\Sigma}) \cong H^0(\Theta_{P_\Sigma}(-\log D)) \oplus \bigoplus_{\sigma\in \Sigma(1)} H^0({\cal O}(D_\sigma)|_{D_\sigma}) 
\]
with the first summand generated by $v_n$ for $n\in N$ and the second summand having a basis $\{mv_{e_\sigma}\}_{m\in \Delta\backslash \Delta^0}$, where $\langle m,e_\sigma\rangle =-1$ for a unique $\sigma\in \Sigma(1)$ determined by $m$.\\
\end{co}
{\bf Remark:} The automorphism group of $P_\Sigma$ apriorily has nothing to do with the toric structure on $P_\Sigma$. The effect of the toric structure is to specify a particular maximal torus in the automorphism group.\\\\

\subsection{Slicing theorem}
For the sake of the mirror symmetry, we need to understand complex moduli near the large complex limit. Notice that the large complex limit (corresponds to $0\in \Delta \subset M$) is not invariant under the full automorphism group, and is only invariant under the toric automorphism group (the maximal torus). To explicitly describe the complex moduli near the large complex limit, it is important to find an explicit slice of the space of sections of $L_\Delta$ that contains the large complex limit, is invariant under toric automorphism group and intersects each nearby orbit of the automorphism group at an orbit of the toric automorphism group (the maximal torus). Let $\Delta^0$ denote the set of integral points in the $(r-2)$-skeleton of $\Delta$.\\
\begin{th}
\label{da}
Near the large complex limit, ${\rm Span}_{\bf C}(\Delta^0)$ is a slice of ${\rm Span}_{\bf C}(\Delta)$ that contains the large complex limit, is invariant under toric automorphism group and intersects each nearby orbit of the automorphism group at an orbit of the toric automorphism group (the maximal torus).\\
\end{th}
To understand the construction better, let's analyze the example of cubics in ${\bf CP}^2$ in detail.\\\\  
{\bf Example:} In the case of cubics in ${\bf CP}^2$, the large complex limit corresponds to $z_0z_1z_2$. We want to reduce\\
\begin{eqnarray*}
F(z)&=& a_0z_0^3 + a_1z_1^3 + a_2z_2^3 + \psi z_0z_1z_2\\
&& \hskip 40pt + a_{01}z_1^2z_2 + a_{02}z_1z_2^2\\
&+& a_{10}z_0^2z_2 \hskip 50pt + a_{12}z_0z_2^2\\
&+& a_{20}z_0^2z_1 + a_{21}z_0z_1^2
\end{eqnarray*}
to\\
\[
F_0(z)= a_0z_0^3 + a_1z_1^3 + a_2z_2^3 + \psi z_0z_1z_2
\]
Consider the transformation\\
\[
T_1:
\left(
\begin{array}{l}
z_0\\z_1\\z_2
\end{array}
\right)
\longrightarrow
\left(
\begin{array}{ccc}
1 & -a_{01}/\psi & -a_{02}/\psi \\
- a_{10}/\psi &1& -a_{12}/\psi \\
-a_{20}/\psi & -a_{21}/\psi & 1
\end{array}
\right)
\left(
\begin{array}{l}
z_0\\z_1\\z_2
\end{array}
\right)
\]
Then\\
\[
F_1(z)=F(T_1z) = \psi z_0z_1z_2 + a_0z_0^3 + a_1z_1^3 + a_2z_2^3 +O(1/\psi)
\]
Or write\\
\begin{eqnarray*}
F_1(z)&=& a_0z_0^3 + a_1z_1^3 + a_2z_2^3 + \psi z_0z_1z_2\\
&& \hskip 40pt + a_{01}z_1^2z_2 + a_{02}z_1z_2^2\\
&+& a_{10}z_0^2z_2 \hskip 50pt + a_{12}z_0z_2^2\\
&+& a_{20}z_0^2z_1 + a_{21}z_0z_1^2
\end{eqnarray*}
with $a_{jk}=O(1/\psi)$. (Here for simplicity, we are abusing the notation a little bit by using the same notation for the coefficients of $F_1$ as for the coefficients of $F$.) Repeat this process inductively, we will get $T_2,T_3,\cdots$ such that\\
\begin{eqnarray*}
F_l(z)=F_{l-1}(T_lz)&=& a_0z_0^3 + a_1z_1^3 + a_2z_2^3 + \psi z_0z_1z_2\\
&& \hskip 40pt + a_{01}z_1^2z_2 + a_{02}z_1z_2^2\\
&+& a_{10}z_0^2z_2 \hskip 50pt + a_{12}z_0z_2^2\\
&+& a_{20}z_0^2z_1 + a_{21}z_0z_1^2
\end{eqnarray*}
with $a_{jk}=O(1/\psi^l)$. Let $T= T_1T_2\cdots$, we get\\
\[
F_0(z) = F(Tz) = a_0z_0^3 + a_1z_1^3 + a_2z_2^3 + \psi z_0z_1z_2
\]
To prove theorem \ref{da} in general, we need to understand the action of $mv_{e_\sigma}$ on sections of the anti-canonical bundle $L_\Delta$. Section of $L_\Delta$, in general, can be expressed as $m(v_{n_1}\wedge \cdots \wedge v_{n_r})$, for $m\in \Delta$, where $n_1,\cdots,n_r$ is a basis of $N$. We have\\
\begin{lm}
\label{df}
\[
{\cal L}_{m_1v_{n_1}} m_2(v_{n_1}\wedge \cdots \wedge v_{n_r})=\langle m_2-m_1,n_1\rangle (m_1+m_2)(v_{n_1}\wedge \cdots \wedge v_{n_r})
\]\\
\end{lm}
In particular, for $mv_{e_\sigma}$, $m(\not=0)\in\Delta\backslash\Delta^0$, $\langle m,e_\sigma\rangle=-1$,\\
\[
{\cal L}_{mv_{e_\sigma}} m_1(v_{n_1}\wedge \cdots \wedge v_{n_r})=(\langle m_1,e_\sigma\rangle+1) (m+m_1)(v_{n_1}\wedge \cdots \wedge v_{n_r})
\]\\
\[
{\cal L}_{mv_{e_\sigma}} (v_{n_1}\wedge \cdots \wedge v_{n_r})= m(v_{n_1}\wedge \cdots \wedge v_{n_r})
\]\\
To simplify the notation, let $s_m = m(v_{n_1}\wedge \cdots \wedge v_{n_r})$,  for $m\in M$ and $v_m = mv_{e_\sigma}$ for $m(\not=0)\in\Delta\backslash\Delta^0$, $\langle m,e_\sigma\rangle=-1$. Then the above formula can be expressed as\\
\begin{co}
\label{dd}
For $m(\not=0)\in\Delta\backslash\Delta^0$, $m_1\in M$, 
\[
{\cal L}_{v_m} s_0=  s_m,\ \ \ {\cal L}_{v_m} s_{m_1}= (\langle m_1,e_\sigma\rangle+1) s_{m+m_1}
\]
where $\sigma\in \Sigma(1)$ is determined by $\langle m,e_\sigma\rangle=-1$.\\
\end{co}
{\bf Proof of theorem \ref{da}:} Let\\
\[
s = \sum_{m(\not= 0)\in \Delta} a_ms_m + \psi s_0,
\]\\
Then the theorem asserts that when $\psi$ is large, automorphism of $P_\Sigma$ can reduce $s$ to the following standard form\\
\[
s^0=\sum_{m \in \Delta^0} a_ms_m + \psi s_0,
\]\\
Consider the automorphism of $P_\Sigma$\\
\[
T_1 = \exp (-\sum_{m(\not=0)\in\Delta\backslash\Delta^0} a_m v_m)
\]\\
Then by using corollary \ref{dd}, we have\\
\begin{eqnarray*}
s^1=T_1(s) &=& \psi s_0 + \sum_{m(\not=0)\in \Delta} a_ms_m - \sum_{m(\not=0)\in\Delta\backslash\Delta^0} a_m {\cal L}_{v_m}s_0+ O(\frac{1}{\psi})\\
&=& \psi s_0 + \sum_{m \in \Delta^0} a_ms_m + \sum_{m(\not=0)\in\Delta\backslash\Delta^0} a_ms_m - \sum_{m(\not=0)\in\Delta\backslash\Delta^0} a_ms_m + O(\frac{1}{\psi})\\
&=& \psi s_0 + \sum_{m \in \Delta^0} a_m s_m + O(\frac{1}{\psi})
\end{eqnarray*}
Adjust $\psi$ and $a_m$ accordingly, above equation can also be expressed as\\
\[
s^1= \sum_{m(\not=0)\in \Delta} a_ms_m + \psi s_0,
\]\\
where $a_m=O(1/\psi)$ for $m(\not=0)\in \Delta\backslash \Delta^0$. Repeat this process inductively, we will get $T_2, T_3,\cdots$ such that\\
\[
s^l = T_l(s^{l-1})= \sum_{m(\not=0)\in \Delta} a_ms_m + \psi s_0,
\]\\
with $a_m=O(1/\psi^l)$ for $m(\not=0)\in \Delta\backslash \Delta^0$. Let $T= T_1T_2\cdots$, we get\\
\[
s^0 = T(s) = \sum_{m \in \Delta^0} a_ms_m + \psi s_0
\]\\
is in standard form (belong to our slice).
\begin{flushright} $\Box$ \end{flushright}

\se{Complex moduli of Calabi-Yau hypersurface and the monomial-divisor mirror map}
According to theorem \ref{da}, the complex moduli near the large complex limit can be characterized by the space of sections\\
\[
s = s_0 + \sum_{m \in \Delta^0} a_ms_m
\]\\
of $L_\Delta$ for $a_m$ small modulo the action of $T=U_0=N_{\bf C}/N$. Let\\
\[
\tilde{M}_0 = \left\{a^I= \prod_{m\in \Delta^0}a_m^{i_m}\ \left|\ I=(i_m)_{m\in \Delta^0}\in {\bf Z}^{\Delta^0}\right.\right\}\cong {\bf Z}^{\Delta^0},
\]\\
then its dual $\tilde{N}_0$ is naturally the space of all the weights\\
\[
\tilde{N}_0 = \{ w=(w_m)_{m\in \Delta^0}\in {\bf Z}^{\Delta^0} \} \cong {\bf Z}^{\Delta^0}.
\]\\
By lemma \ref{df}, we have\\
\[
{\cal L}_{v_n} s_m = \langle m,n\rangle s_m
\]\\
$n\rightarrow \langle m,n\rangle$ define an embedding $N\hookrightarrow \tilde{N}_0$. Let $W\in N_0$ be the image of this embedding. Then $\tilde{M}=W^\perp$ is the space of $T$ invariant functions on complex torus $\tilde{T}_0={\rm Spec}({\bf C}[\tilde{M}_0]) \cong ({\bf C^*})^{\Delta^0}$. $\tilde{N} = \tilde{N}_0/W$ is dual to $\tilde{M}$. For a convex cone $\sigma_0\subset \tilde{N}_0$, let $\sigma^\vee = \sigma_0^\vee\cap \tilde{M}$, then $\sigma = (\sigma^\vee)^\vee$ is naturally the projection of $\sigma_0$ to $\tilde{N}$.\\\\
Let\\
\[
\sigma_0 = \{ w=(w_m)_{m\in \Delta^0}\in {\bf Z}^{\Delta^0}_{\geq 0}\}\subset \tilde{N}_0.
\]\\
Then the fan $\tilde{\Sigma}_0$ consists of subcones of $\sigma_0$ corresponds to the affine toric variety ${\bf C}^{\Delta^0}$. Subcones $\sigma_{S,0}$ of $\sigma_0$ is indexed by subsets $S\subset \Delta^0$.\\
\[
\sigma_{S,0} = \{w\in \sigma_0 | w_m=0,\ {\rm for} \ m\in S\}
\]\\
Let $\sigma_S$, $\sigma$ be the cones that are the projection to $N$ of the cones $\sigma_{S,0}$, $\sigma_0$ in $\tilde{\Sigma}_0$. Notice that $\tilde{N}_0$ can be interpreted as the space of restriction to $\Delta^0$ of the piecewise linear functions on $M$, which we encountered, when we discussed movable cone of a Calabi-Yau hypersurface. This can be easyly seen by defining $p_w^0(m)=-w_m$. In this way, $\sigma_0$ corresponds to piecewise linear functions that are non-positive, namely, convex at origion. $\sigma_{S,0}$ corresponds to piecewise linear functions that are non-positive and vanish at $S$, namely, convex at the cone generated by $S$ (if there are no more other vertices in this cone). $W$ as image of $N$ corresponds to the space of linear function on $M$. $\tilde{N}$ can be viewed as space of restriction to $\Delta^0$ of the piecewise linear functions modulo linear functions on $M$. $\sigma_S$, $\sigma$ also have corresponding convexity interpretation modulo linear functions.\\\\
We are interested in the complex moduli, which is, near the large complex limit, by our slicing theorem, the quotient of ${\bf C}^{\Delta^0}$ by the action of $T$. But the set theoritical quotient is usually not seperable, let along an algebraic variety. To get a seperable algebraic variety as quotient, we have to take the quotient in geomeric invariant theory sense. Namely, we may need to throw away some of the subtorus and collapse some non-generic orbits together when we do the quotient. Then there is always the problem that what to throw away and what to collapse. There are many ways to do it. For our purpose, we want to find a ``canonical'' quotient that suits our need for mirror symmetry.\\\\
There are some general guideline for the choice, we usually want to keep the good generic orbits, in particular, the stablizer is better trivial. A sequence of generic orbits at limit, could converge to a union of several orbits. If some of these limiting orbits are to be keeped, this orbits should be collapsed to one point in the quotient.\\\\
{\bf Example:} ${\bf C^*}$ acts on ${\bf C^2}$ by\\
\[
t\circ (z_1,z_2) = (tz_1,t^{-1}z_2).
\]\\
the generic orbits are $O_c=\{z_1z_2=c\}$. When $c$ approach zero, $O_c$ approach the union of three orbits\\
\[
\{z_1=z_2=0\}\cup \{z_1=0,z_2\not=0\} \cup \{z_2=0,z_1\not=0\}
\]\\
The first orbit has the stablizer ${\bf C^*}$ need to be throw away, the last two orbits have trivial stablizer and can be keeped, but have to be collapsed to one point in the quotient. In this way we get the quotient to be ${\bf C}$ parametrized by $c$.\\\\
According to lemma \ref{bc}, the stablizer of the action of $T$ on $T_{\sigma_{S,0}}$ is corresponding to $W\cap {\rm Span}_{\bf C} (\sigma_{S,0})$. The stablizer is trivial if and only if $W$ projects injectively into $\tilde{N}_{\sigma_{S,0},0}= \tilde{N}_0/{\rm Span}_{\bf C} (\sigma_{S,0})$. In another word,\\
\begin{lm}
The stablizer of the action of $T$ on $T_{\sigma_{S,0}}$ is trivial if and only if the cone generated by $S\in \Delta^0$ is of top dimension in $M$.\\
\end{lm}
In general, we are intersted in the kind of toric varieties whose fan $\Sigma$ is determined by the top dimensional cones $\Sigma(r)$, namely each cone in $\Sigma$ is a subcone of a top dimensional cone in $\Sigma$. Assume our quotient is a toric variety of this type. Then we only need to determine top dimensional cones in its fan $\tilde{\Sigma}$.\\
\begin{lm}
$\sigma_{S,0}$ is of top dimensional in $\tilde{N}$ and $T$ acts freely on $T_{\sigma_{S,0}}$ if and only if\\
\[
\tilde{N}_0 = W \oplus {\rm Span}_{\bf C} (\sigma_{S,0})
\]\\
or equivalently, $S$ contain exactly $r$ linear independent elements in $M$.\\
\end{lm}
For our purpose, we want our construction to be compatible with the Newton polyhedron $\Delta$, namely, we require $S$ belongs to one of the $(r-1)$-dimensional face of $\Delta$. 
This is the way complex moduli of Calabi-Yau hypersurface in $P_\Sigma$ and \k moduli of Calabi-Yau hypersurface in $P_{\Sigma^\vee}$ are related. Namely, the complex moduli of Calabi-Yau hypersurface in $P_\Sigma$ near the large complex limit is naturally a toric variety, its fan can be naturally identified with the fan made up with the varieous \k cones of different birational models of Calabi-Yau hypersurface in $P_{\Sigma^\vee}$.\\\\
The mirror symmetry actually requires something more, a precise identification of complex moduli $P_{\tilde{\Sigma}}$ near the large complex limit and the complexified \k moduli $(\tilde{N}\otimes_{\bf Z}{\bf R} + iK^0_{\Delta^\vee})/\tilde{N}$ near the large radius limit.\\\\
In general, recall that for $n\in N_{\bf C} = N\otimes_{\bf Z}{\bf C}$, $m\in M$, we have the Lie algebra action $v_n(m)= \langle m,n\rangle m$. By the expoential map, we have the Lie group action\\
\[
n(m)=e^{2\pi i\langle m,n\rangle}m
\]\\
where we assume that $n\in T= (N\otimes_{\bf Z}{\bf C})/N$. For any $\sigma\in \Sigma$, we can define $\sigma_{\bf C}= N_{\bf R} + i\sigma_{\bf R} \subset N_{\bf C}$, where $N_{\bf R}= N\otimes_{\bf Z}{\bf R}$ and $\sigma_{\bf R}$ is the real cone in $N_{\bf R}$ generated by $\sigma$. $\sigma^\vee \subset M$ are the monomials that can be extended to the affine toric variety $U_\sigma = {\rm Spec}({\bf C}[\sigma^\vee])$. $\sigma^\vee\backslash \sigma^\perp \subset M$ are the functions on $U_\sigma$ that vanish along $T_\sigma\hookrightarrow U_\sigma$. For $n=n_1 +in_2\in \sigma^0_{\bf C}$, $m\in \sigma^\vee\backslash \sigma^\perp$, we have $n_2\in \sigma^0_{\bf R}$ and $\langle m,n_2\rangle>0$. For $t\in {\bf R}_+$, we have\\
\[
(tn)(m)= e^{2\pi i\langle m,tn\rangle}m = e^{2\pi i\langle m,tn_1\rangle} e^{-2\pi t\langle m,n_2\rangle}m 
\]\\
clearly\\
\[
\lim_{t\rightarrow +\infty} (tn)(m)= 0
\]\\
This implies that for any $n= \in \sigma^0_{\bf C}$ and $x\in T=U_0$, the ray $tn(x)$ approaches a point in $T_\sigma$ when $t\in {\bf R}_+$ approaches $+\infty$. The limit is actually indipendent of $n \in \sigma^0_{\bf C}$.\\
\[
\lim_{t\rightarrow +\infty} (tn)(x)= x_\sigma
\]\\
where $x_\sigma$ is the image of $x$ under the projection $\pi_\sigma : T= {\rm Spec}({\bf C}[M]) \rightarrow T_\sigma = {\rm Spec}({\bf C}[\sigma^\perp])$ that is determined by the embeding $\sigma^\perp \hookrightarrow M$. Actually $n(x)\in \pi^{-1}_\sigma(x_\sigma)$ for $n\in \sigma_{\bf C}\subset \bar{\sigma}_{\bf C}$.\\\\ 
Apply this construction to our situation, we get the monomial divisor map\\
\[
(\tilde{N}\otimes_{\bf Z}{\bf R} + iK^0_{\Delta^\vee})/\tilde{N}\rightarrow P_{\tilde{\Sigma}}
\]\\
which is an approximation of the actual mirror map.\\\\

{\bf Question:} Fan $\tilde{\Sigma}$ determine a sense of being near the large complex limit. What is precise characterization of being near the large complex limit determined by the fan $\tilde{\Sigma}$?\\\\
  
For $\sigma\in \Sigma$, $\sigma_{\bf C}$ is a subsemigroup of the complex torus $T\cong U_0$. For any $x\in U_0$, $O_x= \sigma_{\bf C}(x)$ is an orbit of $\sigma_{\bf C}$. Two orbits $O_x$, $O_y$ of $\sigma_{\bf C}$ are called equivalent, if $O_x\cap O_y$ is non-empty. We denote the equivalent class by $[O_x]$. Use these equivalence classes of orbits of semigroups, we have another interpretation of $P_\Sigma$.\\
\[
P_\Sigma \cong \{[O]|O\subset U_0\ {\rm is\ an\ orbit\ of\ } \sigma_{\bf C}\ {\rm for\ some\ } \sigma\in \Sigma\}
\]\\
From this point of view, for $w=(w_m)_{m\in\Delta^0} \in K^0_{\Delta^\vee}$, and\\
\[
s = s_0 + \sum_{m\in \Delta^0} a_ms_m
\]\\ 
\[
t^w\circ s = s_0 + \sum_{m\in \Delta^0} (t^{w_m}a_m)s_m
\]\\
should be near the large complex limit if $t\in {\bf R}_+$ is sufficiently small. In another word, {\bf for $a=(a_m)_{m\in\Delta^0}$ small,\\
\[
s = s_0 + \sum_{m\in \Delta^0} a_ms_m
\]\\ 
is near the large complex limit $s_0$ if $(-\log a_m)_{m\in\Delta^0}$ is in $(K^0_{\Delta^\vee})_{\bf C}$, or equivalently, $(-\log |a_m|)_{m\in\Delta^0}$ is in $K^0_{\Delta^\vee}$.}\\\\

\se{Lagrangian torus fibration of Calabi-Yau hypersurface}
In \cite{lag1}, \cite{lag2}, we constructed Lagrangian fiberation for quintic Fermat type Calabi-Yau familly $\{X_\psi\}$ in ${\bf P^4}$ by gradient flow method. In \cite{lag3}, we deal with the case of generic quintic. In generic quintic case, the difficulty is not gradient flow. In very general situation, based on the technique result in \cite{lag2}, the gradient flow will deform the canonical Lagrangian torus fibration of the large complex limit to form Lagrangian torus fibration on general quintic Calabi-Yau. But the singular locus of the fibration so constructed is usually of codimension 1, while the SYZ conjecture expects that the singular locus to be of codimension 2. In the case of Fermat type quintic case, due to the explicit nature of the singular locus, it is easy to observe that the singular locus is naturally a fattening of a 1-dimensional graph $\Gamma$. Use some explicit computation and symplectic geometry technique, in \cite{lag2} we were able to squize the singular locus to construct a Lagrangian torus fibration with codimension 2 singular locus $\Gamma$. In the generic quintic case, the singular locus is usually very complicated and in general does not necessaryly resemble a fattening of a 1-dimensional graph. However, as observed in \cite{lag3}, when the quintic Calabi-Yau is near the large complex limit in a suitable sense, miraculously, the singular locus will resemble a fattening of a 1-dimensional graph $\Gamma$. In \cite{N}, we discussed this phenomenon in great detail and used similar technique as in \cite{lag2} to finally achieve codimension 2 singular locus for Lagrangian torus fibration of generic Calabi-Yau quintic near the large complex limit in \cite{lag3}.\\\\
In this section, we will discuss how to generalize the results in \cite{lag1}, \cite{lag2}, \cite{lag3}, to construct Lagrangian torus fibration with codimension 2 singular locus for generic Calabi-Yau hypersurface in toric variety. Although Calabi-Yau hypersurface in general toric variety sounds much more general than quintic in ${\bf P^4}$. Actually there are no essentially new technical difficulties for the gradient flow. Generalization is a very natural one.\\\\
\subsection{Lagrangian fibration with codimension 1 singular locus}
Before getting into the detail, let's first fix some notation. Assume that $\Delta$ is a reflexive polyhedron, $p_w$ ($p_v$) is a convex piecewise linear function on $M$ ($N$) with respect to $\Delta^0$ ($(\Delta^\vee)^0$). $p_w(m)=-w_m$ ($p_v(n)=-v_n$) for $m\in\Delta^0$ ($n\in(\Delta^\vee)^0$). $p_w$ ($p_v$) naturally determines a fan $\Sigma^w$ ($\Sigma^v$) for $M$ ($N$) that is compatible with $\Delta^\vee$ ($\Delta$). $p_w$ ($p_v$) also naturally determines a polyhedron $\Delta_w\subset N$ ($\Delta_v\subset M$) consists of $n\in N$ ($m\in M$) that as a linear function on $M$ ($N$) is greater or equal to $p_w$ ($p_v$). In particular, the fan $\Sigma$ ($\Sigma^\vee$) of the anti-canonical model $P_{\Sigma}$ ($P_{\Sigma^\vee}$) is determined by $p_\Delta$ ($p_{\Delta^\vee}$), $p_\Delta(n)=-1$ ($p_{\Delta^\vee}(m)=-1$) for $n\in{\Delta^\vee}\backslash\{0\}$ ($m\in\Delta\backslash\{0\}$), for which the corresponding polyhedron is $\Delta$ ($\Delta^\vee$).\\\\
The fan $\Sigma^v$ ($\Sigma^w$) is a crepant subdivision of the fan $\Sigma$ ($\Sigma^\vee$). The corresponding toric variety $P_{\Sigma^v}$ ($P_{\Sigma^w}$) is a crepant resolution of the anti-canonical model $P_{\Sigma}$ ($P_{\Sigma^\vee}$), with the natural map\\
\[
\pi_v: P_{\Sigma^v} \rightarrow P_{\Sigma}
\]
\[
\pi_w: P_{\Sigma^w} \rightarrow P_{\Sigma^\vee}
\]\\
Before getting into the gradient flow construction, it is important to construct the suitable $T_{\bf R}$-invariant \k metric $g$ on $P_{\Sigma^v}$. According to philosophy of SYZ conjecture, our purpose is not merely getting Lagrangian fibration for $X_s$. Instead, we need to construct Lagrangian fibration for $X_s$ that is depnding suitablly on any specific complex structure and \k structure on $X_s$. Our construction of $T_{\bf R}$-invariant \k metric $g$ on $P_\Sigma$ should reflect this philosophy.\\\\
In the case of quintic Calabi-Yau that we dealt with in our previous papers, we essentially used the Fubini-Study metric. The reason is that the \k moduli for quintic Calabi-Yau is of dimension 1. The restriction of the Fubini-Study metric to quintic Calabi-Yau is enough to capture the \k moduli of quintic. Even in that case, to make the singular locus better understood, it is a good idea to modify the the Fubini-Study metric in the same \k class according to the complex structure of the particular quintic Calabi-Yau. It is helpful to review this construction for the quintic Calabi-Yau here to help us to understand our construction of $T_{\bf R}$-invariant \k metric here in the general toric case.\\\\
Consider a generic quintic defined by\\
\[
p(z) = \sum_{|m|=5}a_mz^m=0
\]\\
in ${\bf P^4}$. Let $\Delta$ denote the Newton polygon for quintic polynomials. Then we can introduce the \k potential\\
\[
h_a = \log(|z|_a^2) = \log \left( \sum_{m\in \Delta} |a_mz^m|^2\right).
\]\\
Let $z^m = |z^m|e^{i\theta_m}$, then\\
\[
\omega = \partial \bar{\partial} h_a = i\sum_{m\in \Delta} d\theta_m \wedge d\left(\frac{|a_mz^m|^2}{|z|_a^2}\right)
\]\\
define a $T_{\bf R}$-invariant \k form on ${\bf P^4}$ representing the \k class corresponding to ${\cal O}(5)$. Notice that $\Delta$ is the Newton polygon of ${\cal O}(5)$ on ${\bf P^4}$.\\
\begin{lm}
The moment map is\\
\[
F_a(x) = \sum_{m\in \Delta} \frac{|a_mz^m|^2}{|z|_a^2}m,
\]\\
which maps ${\bf P^4}$ to $\Delta$.\\
\end{lm}
Near the large complex limit, under this moment map, the singular locus behave very nicely. This enable us to construct Lagrangian torus fibration with codimension 2 singular locus for generic quintic Calabi-Yau in \cite{lag3}.\\\\
For the case of general Calabi-Yau hypersurface in toric variety, $p_w$ ($p_v$) as a convex piecewise linear function on $M$ ($N$) with respect to $\Delta^*$ (${\Delta^\vee}^*$) can be naturally understood as a \k class for Calabi-Yau hypersurface in $P_{\Sigma^w}$ ($P_{\Sigma^v}$). Recall that under the monomial-divisor mirror map, the complex moduli of Calabi-Yau hypersurface in $P_{\Sigma^v}$ is naturally identified with complexified \k cone of Calabi-Yau hypersurface in $P_{\Sigma^w}$. Therefore $p_w$ also characterize complex moduli of Calabi-Yau hypersurface in $P_{\Sigma^v}$. More precisely, for the Calabi-Yau hypersurface $X_s$ defined by\\
\[
s = \sum_{m(\not= 0)\in \Delta} a_ms_m + \psi s_0 =0,
\]\\
\[
p_w(m) = -w_m = \log |a_m/\psi|.
\]\\
The pair $(p_w,p_v)$ together characterize the whole moduli (including both the complex and the \k moduli) of a Calabi-Yau hypersurface in $P_{\Sigma^v}$. For the purpose of our gradient flow, we need to construct a $T_{\bf R}$-invariant \k metric $g$ on $P_{\Sigma^v}$ that depends on both $p_v$ and $p_w$. More precisely, we want the \k form of $g$ to correspond to $p_v$. Also to make the singular locus well behaved, we want the \k form to depend on complex moduli $p_w$ in a way similar to the quintic case.\\\\
Now we can define a natural $T$-invariant \k form on $P_{\Sigma^v}$. $m\in M$ naturally define a holomorphic function $s_m$ on $T=N_{\bf C}/N$, $|s_m|$ is a $T_{\bf R}$-invariant function on $T$. For general $m\in M_{\bf R}$, $s_m$ will not make sense as function, but we always can define $|s_m|$ as a $T_{\bf R}$-invariant function on $T$. Let\\
\[
h_{w,v}=\log(|s|_{w,v}^2) = \log\left(\int_{m\in\Delta_v}|s_m|_w^2\right),
\]\\
where $|s_m|_w^2=|e^{-w_m}s_m|^2$. 
\[
\omega = \partial \bar{\partial} h_{w,v} = i\int_{m\in\Delta_v} d\langle m,\theta\rangle \wedge d\left(\frac{|s_m|_w^2}{|s|_{w,v}^2}\right)
\]\\
define a $T$-invariant \k form on $P_{\Sigma^v}$ representing the \k class corresponding to $p_v$.\\
\begin{lm}
The moment map is\\
\[
F_{w,v}(x) = \int_{m\in\Delta_v} \frac{|s_m|_w^2}{|s|_{w,v}^2}m,
\]\\
which maps $P_{\Sigma^v}$ to $\Delta_v$.\\
\end{lm}
By this map, $T_{\bf R}$-invariant functions $h_{w,v}$, $\frac{|s_m|_w^2}{|s|_{w,v}^2}$ can all be viewed as function on $\Delta_v$. We have\\
\begin{lm}
$\rho_m = \frac{|s_m|_w^2}{|s|_{w,v}^2}$ as a function on $\Delta_v$ archieves maximum exactly at $m\in \Delta_v$.\\
\end{lm}
{\bf Proof:} By\\
\[
x_k\frac{\partial |s_m|_w^2}{\partial x_k} = \langle m,n_k\rangle |s_m|_w^2
\]\\
$\rho_m = \frac{|s_m|_w^2}{|s|_{w,v}^2}$ archieves maximal implies\\
\[
\sum_{k=1}^r x_k\frac{\partial \rho_m}{\partial x_k}m_k = \frac{|s_m|_w^2}{|s|_{w,v}^2}\int_{m'\in\Delta}(\langle m,n_k\rangle - \langle m',n_k\rangle)\frac{|s_{m'}|_w^2}{|s|_{w,v}^2}m_k
\]
\[
= \frac{|s_m|_w^2}{|s|_{w,v}^2}\int_{m'\in\Delta}(m - m')\frac{|s_{m'}|_w^2}{|s|_{w,v}^2} = \frac{|s_m|_w^2}{|s|_{w,v}^2}(m - F_{w,v}(x))=0
\]\\
Therefore\\
\[
F_{w,v}(x) =m
\]\\
when $\rho_m = \frac{|s_m|_w^2}{|s|_{w,v}^2}$ archieves maximal.
\begin{flushright} $\Box$ \end{flushright}
The moment maps naturally induce the following diagram\\
\begin{center}
\setlength{\unitlength}{1.2pt}
\begin{picture}(100,50)(0,0)
\put(30,40){\makebox(0,0){$P_{\Sigma^v}$}}
\put(80,40){\makebox(0,0){$P_{\Sigma}$}}
\put(30,10){\makebox(0,0){$\Delta_v$}}
\put(80,10){\makebox(0,0){$\Delta$}}
\put(42,40){\vector(1,0){26}}
\put(28,32){\vector(0,-1){15}}
\put(42,10){\vector(1,0){26}}
\put(80,32){\vector(0,-1){15}}
\put(55,45){\makebox(0,0){\footnotesize{$\pi_v$}}}
\put(55,15){\makebox(0,0){\footnotesize{$\hat{\pi}_v$}}}
\put(38,25){\makebox(0,0){\footnotesize{$F_{w,v}$}}}
\put(85,25){\makebox(0,0){\footnotesize{$F$}}}
\end{picture}
\end{center}

Recall that $L_\Delta$ is the anti-canonical bundle of $P_\Delta$. $s_m$ for $m\in \Delta$ give us the set of $T$-invariant sections of $L_\Delta$. $\pi_v^*L_\Delta$ is the anti-canonical bundle of $P_{\Sigma^v}$. Consider a generic section of $L_\Delta$\\
\[
s = \sum_{m(\not= 0)\in \Delta} a_ms_m + \psi s_0,
\]\\
$s$ can also be thought of as a section of $\pi_v^*L_\Delta$. Zero set $X_s$ of $s$ as a section of $\pi_v^*L_\Delta$ is a Calabi-Yau hypersurface in $P_{\Sigma^v}$. When we vary $s$, we get the family of Calabi-Yau hypersurfaces corresponds to sections of $L_\Delta$. When $\psi$ approach $\infty$, the corresponding Calabi-Yau hypersurface $X_s$ approach the ``large complex limit" $X_\infty$ defined by $s_0=0$. The moment map $F_{w,v}$ naturally give us the Lagrangian fibration\\
\[
F_{w,v}: X_\infty \rightarrow \partial \Delta_v
\]\\
of the large complex limit $X_\infty$.\\\\
Consider the meromorphic function\\
\[
q(z)=\frac{s_0}{s-\psi s_0}
\]\\
defined on $P_\Sigma$. Let $\omega$ denote the \k form of a $T_{\bf R}$-invariant \k metric $g$ on $P_\Sigma$, and $\nabla f$ denote the gradient vector field of real function $f=Re(q)$ with respect to the \k metric $g$. Clearly, $X_\infty$ is invariant under the action of $T_{\bf R}$. The fibres of the natural Lagrangian fibration of $X_\infty$ are exactly the orbits of $T_{\bf R}$-action. From our previous discussion, we know that the flow of $V=\frac{\nabla f}{|\nabla f|^2}$ moves $X_{\infty}$ to $X_s$ symplectically, therefore gives rise to Lagrangian fiberation $F_s$ of $X_s$ over $\Delta_v$.\\\\
From our previous experience with Fermat type quintic Calabi-Yau and generic quintic Calabi-Yau, we know that the Lagrangian fibration we get by this way is not a $C^{\infty}$ Lagrangian fibration and usually with wrong type of singular fibers.\\\\
In the Fermat type quintic case, due to the highly symmetric nature of Fermat type quintic Calabi-Yau, it was relatively easy to figure out what is the expected $C^{\infty}$ Lagrangian (even special Lagrangian) fibration structure and corresponding singular fibres. Also in that case, it is easy to see explicitly that the singular locus of the Lagrangian fibration we constructed by gradient flow is a fattened version of the singular locus of the expected Lagrangian fibration. Therefore, even without perturbing to the $C^{\infty}$ Lagrangian fibration, we can already compute the expected monodromy.\\\\  
In the case of general Calabi-Yau hypersurfaces in toric variety, just like the case of generic quintic Calabi-Yau discussed in \cite{lag3}, to carry out similar construction we just described for the Fermat type quintic, there are several obstacles. One of the major difficulty to generalize the discussion to the case of general Calabi-Yau hypersurfaces in toric variety is that for general Calabi-Yau hypersurfaces in toric variety the singular locus of the Lagrangian fibration constructed from deforming the standard Lagrangian fibration of $X_{\infty}$ via the flow of $V$ can be fairly arbitrary and does not necessarily resemble the fattening of any ``expected singular locus''. Worst of all, in the case of general Calabi-Yau hypersurfaces in toric variety, there are no obvious guess what the ``expected'' $C^{\infty}$ Lagrangian (special Lagrangian) fiberation should be. One clearly expect the expected singular locus to be some graph in $\partial \Delta \cong S^3$. But it seems to take some miracle (at least to me when I first dream about it) for a general singular set $C$ (which is an algebraic curve) to project to the singular locus $\tilde{\Gamma} = F(C)$ that resemble a fattened graph.\\\\
Interestingly, miracle happens here! It largely relies on better understanding of what it means to be {\bf near the large complex limit}. Philosophically speaking, it is commonly believed that the large complex limit corresponds to classical situation (in comparison to quantum situation). The closer one gets to the large complex limit, the smaller the quantum effects will be. Therefore we may explain the possible chaotic behavior of the singular locus $\tilde{\Gamma} = F(C)$ by quantum effect. To make the singular locus behave nicely to resemble a one-dimensional graph, we simply need to get really close to the large complex limit $X_{\infty}$. Therefore the qestion remain is: how to get close to the large complex limit?\\\\ 
Before getting into the detail, let's recall facts of our gradient flow. One way to understand $V$ is to realize that $dq$ is a meromorphic section of $N^*_X$, $(dq)^{-1}$ is a natural meromorphic section of $N_X$. Notice the exact sequence\\
\[
0 \rightarrow T_X \rightarrow T_{P_\Sigma}|_X \rightarrow N_X \rightarrow 0
\]\\
With respect to the \k metric $g$ on $T_{P_\Sigma}$, the exact sequence has a natural (non-holomorphic) splitting\\
\[
T_{P_\Sigma}|_X =T_X \oplus N_X
\]\\
$V$ is just real part of the natural lift of $(dq)^{-1}$ via this splitting. $V$ is singular exactly when $(dq)^{-1}$ is singular or equivalently, when $dq=0$, which corresponds to singular part of $X_{\infty}$. On the other hand, the union of the Lagrangian 3-torus non-singular fibers of $X_{\infty}$ is exactly the smooth part of $X_{\infty}$. Since $V$ is non-singular on the smooth part of $X_{\infty}$, all these 3-torus fibers will be carried to $X_s$ nicely by the flow of $V$.\\\\
Now we will try to understand how the gradient flow of $f$ behaves at singularities of $X_{\infty}$. As we know, $X_{\infty}$ has only normal crossing singularities. ${\rm Sing}(X_{\infty})$ is a union of complex torus (orbits of $T$) of dimension less or equal to $r-2$. According to general results in \cite{lag2} about the gradient vector field, we have\\
\begin{prop}
\label{fa}
Under the flow of $V$, points in $X_{\infty}\cap X_s$ will be fixed. Points in ${\rm Sing}(X_{\infty})\backslash X_s$ that is located in a $l$-dimensional complex torus in ${\rm Sing}(X_{\infty})$ will bubble off to a $(r-l-1)$-dimensional real torus in $X_s$.\\
\end{prop}
From this general structure result on gradient flow, it is easy to pin down the topological singular set $C$ of the fibration $F_s$.\\
\begin{co}
The topological singular set of the Lagrangian torus fibration $F_s$ is $C = {\rm Sing}(X_{\infty})\cap X_s$, and the corresponding singular locus is $\tilde{\Gamma} = F_s(C)$.\\
\end{co}
Recall that $T_{\bf R}$-orbits form the fibres of the Lagrangian fibration $F_{w,v}: X_\infty \rightarrow \partial\Delta_v$. By proposition \ref{fa}, a $T_{\bf R}$-orbit in ${\rm Sing}(X_{\infty})$ that does not intersect $C$ will deform to a Lagrangian 3-torus in $X_s$ under the flow of $V$. $T_{\bf R}$-orbits intersecting $C$ will deform to singular fibres of the Lagrangian fibration $F_s: X_s \rightarrow \partial\Delta_v$. The actural singular fibre type is determined by the way the $T_{\bf R}$-orbit is intersecting $C$. To further discuss singular fibre type, it is important to understand better the singular set $C$ and the corresponding singular locus $\tilde{\Gamma}$.\\\\
In general, for any $\sigma\in \Sigma$, we have the quotient fan $\Sigma_\sigma$, $P_{\Sigma_\sigma}$ is the clossure of $T_\sigma$ in $P_{\Sigma}$. $\Delta_{\sigma^*} = \sigma^*\cap (\Delta\backslash\{0\})$ is a subface of $\Delta$ and is the Newton polyhedron of $T_\sigma$-invariant sections of $L_\Delta|_{P_{\Sigma_\sigma}}$. $s_w$ is restrict to\\
\[
s_{w,\sigma} = \sum_{m\in \Delta_{\sigma^*}} e^{-w_m}s_m
\]\\
on $P_{\Sigma_\sigma}$. Let $C_\sigma = \{s_{w,\sigma}^{-1}(0)\}$.\\\\
The singular set $C$ can be understood as the curve in the 2-skeleton $X^{(2)}_\infty = {\rm Sing}(X_{\infty})$ of $X_\infty$ that is cut out by $s=0$. For 2-skeleton, we have the following diagram\\
\begin{center}
\setlength{\unitlength}{1.3pt}
\begin{picture}(100,50)(0,0)
\put(30,40){\makebox(0,0){$P_{\Sigma^v}^{(2)}$}}
\put(80,40){\makebox(0,0){$P_{\Sigma}^{(2)}$}}
\put(30,10){\makebox(0,0){$\Delta_v^{(2)}$}}
\put(80,10){\makebox(0,0){$\Delta^{(2)}$}}
\put(42,40){\vector(1,0){26}}
\put(27,32){\vector(0,-1){15}}
\put(42,10){\vector(1,0){26}}
\put(77,32){\vector(0,-1){15}}
\put(55,46){\makebox(0,0){\footnotesize{$\pi_v^{(2)}$}}}
\put(55,16){\makebox(0,0){\footnotesize{$\hat{\pi}_v^{(2)}$}}}
\put(36,25){\makebox(0,0){\footnotesize{$F_{w,v}^{(2)}$}}}
\put(85,25){\makebox(0,0){\footnotesize{$F^{(2)}$}}}
\end{picture}
\end{center}
$P_{\Sigma}^{(2)}$ is a union of complex torus of dimension less or equal to 2. Generic section $s$ can avoid verteces (the zero dimensional tori) and will intersect complex torus of dimension 1 and 2. In general\\
\[
\dim \Delta_{\sigma^*} + \dim \sigma = r.
\]\\
Therefore, for $\sigma\in \Sigma(3)$, $\dim_{\bf C}(C_\sigma)=0$, for $\sigma\in \Sigma(2)$, $\dim_{\bf C}(C_\sigma)=1$. Notice that ${\rm Sing}(X_{\infty}) = P_{\Sigma^v}^{(2)}$. The singular set $C$ can be decomposed into two parts $C=C_1\cup C_2$, where\\
\[
C_1 = \bigcup_{\sigma\in \Sigma(3)} (\pi_v^{(2)})^{-1}(C_\sigma)
\]\\ 
\[
C_2 = \bigcup_{\sigma\in \Sigma(2)} (\pi_v^{(2)})^{-1}(C_\sigma)
\]\\
The corresponding singular locus $\tilde{\Gamma}= \tilde{\Gamma}_1 \cup \tilde{\Gamma}_2$. It is easy to see that $\dim \tilde{\Gamma}_1=1$ and $\dim \tilde{\Gamma}_2=2$. Let $\tilde{\Gamma}^0$ denotes the set of verteces of $\tilde{\Gamma}$, $\tilde{\Gamma}_1^1= \tilde{\Gamma}_1\backslash \tilde{\Gamma}^0$, $\tilde{\Gamma}_2^2$ denotes the interior of $\tilde{\Gamma}_2$, and $\tilde{\Gamma}_2^1= \tilde{\Gamma}_2\backslash (\tilde{\Gamma}_2^0 \cup \tilde{\Gamma}_2^2)$. When $X_s$ is near the large complex limit (we will descuss in more detail in a moment), under the moment map, $C_2$ is 2-1 over $\tilde{\Gamma}_2^2$ and 1-1 over $\tilde{\Gamma}_2^1$.\\\\ 
For the fibration $F_{w,v}: X_\infty \rightarrow \partial\Delta_v$, when $X_s$ is near the large complex limit, fibre over a point in $\tilde{\Gamma}^0$ is a circle that intersects $C$ at a point, fibre over a point in $\tilde{\Gamma}_1^1$ is a 2-torus that intersects $C$ at a circle, fibre over a point in $\tilde{\Gamma}_2^1$ is a 2-torus that intersects $C$ at a point, fibre over a point in $\tilde{\Gamma}_2^2$ is a 2-torus that intersects $C$ at 2 points. By proposition \ref{fa}, we can pin down all the singular fibres of the Lagrangian fibration $F_s$.\\
\begin{th}
When $X_s$ is near the large complex limit, flow of $V$ will produce a Lagrangian fibration $F_s: X_s \rightarrow \partial\Delta_v$. There are 5 types of fibers.\\
(i). For $p\in \partial\Delta\backslash \tilde{\Gamma}$, $F^{-1}(p)$ is a Lagrangian 3-torus.\\
(ii). For $p\in \tilde{\Gamma}_2^2$, $F_s^{-1}(p)$ is a Lagrangian 3-torus with $2$ circles collapsed to $2$ singular points.\\ 
(iii). For $p\in \tilde{\Gamma}_2^1$, $F_s^{-1}(p)$ is a Lagrangian 3-torus with $1$ circle collapsed to $1$ singular point.\\
(iv). For $p\in \tilde{\Gamma}^0$, $F_s^{-1}(p)$ is a Lagrangian 3-torus with $1$ 2-torus collapsed to $1$ singular point (the type III singular fibre).\\
(v). For $p\in \tilde{\Gamma}_1^1$, $F_s^{-1}(p)$ is a type I singular fibre (a circle times a rational curve with one node).\\
\end{th}
The Lagrangian torus fibration $F_s: X_s \rightarrow \partial\Delta_v$ we constructed by flow of $V$ is not really the kind of Lagrangian torus fibration that is expected by the SYZ mirror conjecture. SYZ construction generally requires that singular locus of the special Lagrangian fibration to be some one dimensional graph in $S^3$. For our purpose we want to construct $C^\infty$ Lagrangian fibrations with dimension one singular locus that we believe will capture symplectic topological structure of special Lagrangian fibrations in general. Or at least, we want to construct a Lagrangian fibration with dimension one singular locus. To achieve these, in the Fermat type quintic case, the singular locus can easily be seen as a fattening of a one dimensional graph $\Gamma$. In \cite{lag2}, we use symplectic technique to explicitely deform the Lagrangian fibration before we deform acccording the flow to ensure that the singular locus $F(C)=\Gamma$. Then the corresponding Lagrangian fibration of $X_s$ we construct by the flow of $V$ will have dimension one singular locus. In \cite{lag2}, we also discussed the question of how to deform this Lagrangian fibration into a $C^\infty$ Lagrangian fibration.\\\\
In the case of general Calabi-Yau hypersurface in toric variety, we can do exactly the same thing, granted that we can realize the singular locus $F_s(C)$ as a fattening of some graph $\Gamma$ and be able to explicitly construct symplectic deformation, which deform $C$ to satisfy $F_s(C)=\Gamma$. Recall that $C=C_1\cup C_2$, and $\tilde{\Gamma}_1 = F_s(C_1)$ is already of dimension 1. We only need to concentrate on the $C_2$ part. Since $C_2$ is reducible, and each irreducible component is in a 2-dimensional toric subvariety $P_{\Sigma^v_\tau}$ of $P_{\Sigma^v}$ that is over a 2-dimensional toric subvariety $P_{\Sigma_\sigma}$ of $P_{\Sigma}$. according to \cite{lag2}, the problem can be isolated to each $P_{\Sigma^v_\tau}$. More precisely, for $\tau\subset \sigma$, $\tau \in\Sigma^v(2)$, $\sigma \in \Sigma(2)$, $\pi_v: P_{\Sigma^v_\tau} \rightarrow P_{\Sigma_\sigma}$ is a resolution. For $X_s$ generic, the restriction $\pi_v: C_2\cap P_{\Sigma^v_\tau} \rightarrow C_2\cap P_{\Sigma_\sigma}$ is an isomorphism. The problem can further be reduced to discussion over $P_{\Sigma_\sigma} \subset P_{\Sigma}$ as the following.\\\\
{\bf Problem:} Given an integral polyhedron $\Delta \subset M$, one can naturally construct a fan $\Sigma$ and an ample line bundle $L_\Delta$ over the toric variety $P_\Sigma$, such that the Newton polygon of sections of $L_\Delta$ is exactly $\Delta$. Let $F:P_{\Sigma} \rightarrow \Delta$ be the moment map with respect to the class of $L_\Delta$. We need a curve $C_s = \{s^{-1}(0)\}$ in $P_{\Sigma}$, where $s$ is a section of $L_\Delta$, such that $F(C_s)$ is a fattening of some graph $\Gamma$ and we want to explicitly construct symplectic deformation of $P_\Sigma$, which deform $C$ to satisfy $F(C_s)=\Gamma$.\\\\
Clearly, one can not expect that $C_s$ for all section $s$ of $L_\Delta$ will have such a nice property. As we point out earlier, it turns out, when the Calabi-Yau hypersurface in toric variety is generic and close to the large complex limit in a certain sense, the corresponding singular set curves $C_\sigma$ in $P_{\Sigma_\sigma}$ have the properties described in above problem. This kind of curve and more general situation have been discussed intensively in our paper \cite{N}. To describe the result from there, let's first introduce some notation and definations.\\\\
\subsection{Newton polygon and string diagram}
Given an integral polyhedron $\Delta \subset M$, one can naturally construct a fan $\Sigma$ and an ample line bundle $L_\Delta$ over the toric variety $P_\Sigma$, such that the Newton polygon of sections of $L_\Delta$ is exactly $\Delta$. In this context, we assume that $r={\rm rk}(M)=2$.\\\\
Consider a general section of $L_\Delta$\\
\[
s= \sum_{m\in \Delta}a_m s_m
\]\\
Let $|a_m|=e^{-w_m}$, then $w=(w_m)_{m\in \Delta}$ define a function on integral points in $\Delta$ (thinking of $\Delta$ as a real triangle in $M_{\bf R} = M\otimes {\bf R}$).\\\\
\begin{de}
\label{fb}
$w=(w_m)_{m\in \Delta}$ is called convex on $\Delta$, if for any $m' \in \Delta$ there exist an affine function $n$ such that $n(m')= w_{m'}$ and $n(m)\leq w_m$ for $m\in \Delta\backslash \{m'\}$.
\end{de} 
We will always assume that $w$ is convex. With $w$ we can define the moment map\\
\[
F_w(x) = \sum_{m\in\Delta} \frac{|s_m|_w^2}{|x|_w^2}m,
\]\\
which maps $P_\Sigma$ to $\Delta$, where $|s_m|_w^2=|e^{-w_m}s_m|^2$, $|s|_w^2 = \sum_{m\in\Delta}|s_m|_w^2$.\\\\
$\Delta$ can also be thought of as a real triangle in $M_{\bf R} = M\otimes {\bf R}$. Then $w=(w_m)_{m\in \Delta}$ define a function on integral points in $\Delta$. If $w$ is convex, then $w$ can be extended to a piecewise linear convex function on $\Delta$. We will denote the extension also by $w$. Generic $w$ will determine a simplicial decomposition of $\Delta$, with zero simplices being integral points in $\Delta$. In this case, we say the piecewise linear convex function $w$ is compatible with the simplicial decomposition of $\Delta$, and the simplicial decomposition of $\Delta$ is determined by $w$. Consider the baricenter subdivision of this simplicial decomposition of $\Delta$, let $\Gamma_w$ denote the union of simplices in the baricenter subdivision that do not intersect integral points in $\Delta$. Then it is not hard to see that $\Gamma_w$ is an one dimensional graph. And $\Gamma_w$ divide $\Delta$ into regions, with a unique integral point of $\Delta$ located at the center of each of this region. In particular we can think of the regions as parametrized by integral points in $\Delta$.\\\\
Let $C$ denotes the curve defined by $s=0$, the results in \cite{N} give the following.\\
\begin{th}
\label{fc}
For $w=(w_m)_{m\in \Delta}$ convex and positive, and $|w|$ large enough, $F_w(C)$ will be a fattening of graph $\Gamma_w$. There exist a symplectic diffeomorphism that map $C$ to $\tilde{C}$ such that $F_w(\tilde{C})=\Gamma_w$.\\
\end{th}
{\bf Remark:} The moment map $F_w$ is invariant under the real 2-torus action. For any other moment map $F$ that is invariant under the real 2-torus action, we also have that $F(\tilde{C})$ is a 1-dimensional graph.\\\\
{\bf Example:} For the standard simplicial decomposition of the Newton polygon $\Delta_5$ of quintic polyonmials,\\
\begin{center}
\setlength{\unitlength}{1.2pt}
\begin{picture}(200,160)(-30,10)
\put(150,48){\line(-3,-5){18}}
\put(6,48){\line(3,-5){18}}
\put(60,138){\line(1,0){36}}
\put(132,78){\line(-3,-5){36}}
\put(24,78){\line(3,-5){36}}
\put(42,108){\line(1,0){72}}
\put(114,108){\line(-3,-5){54}}
\put(42,108){\line(3,-5){54}}
\put(24,78){\line(1,0){108}}
\put(96,138){\line(-3,-5){72}}
\put(60,138){\line(3,-5){72}}
\put(6,48){\line(1,0){144}}
\put(78,168){\line(-3,-5){90}}
\put(78,168){\line(3,-5){90}}
\put(-12,18){\line(1,0){180}}
\multiput(78,168)(36,0){1}{\circle*{4}}
\multiput(60,138)(36,0){2}{\circle*{4}}
\multiput(42,108)(36,0){3}{\circle*{4}}
\multiput(24,78)(36,0){4}{\circle*{4}}
\multiput(6,48)(36,0){5}{\circle*{4}}
\multiput(-12,18)(36,0){6}{\circle*{4}}
\end{picture}
\end{center}
\begin{center}
\stepcounter{figure}
Figure \thefigure: the standard simplicial decomposition
\end{center}
We have the corresponding $\hat{\Gamma} = F(C)$\\
\begin{center}
\begin{picture}(200,180)(-30,0)
\thicklines
\put(64,158){\line(3,-1){12}}
\put(76,154){\line(1,0){12}}
\put(100,158){\line(-3,-1){12}}

\multiput(0,0)(-18,-30){4}
{\put(64,158){\line(4,-3){15.2}}
\put(79.2,146.6){\line(0,-1){17}}
\put(46,128){\line(3,-1){18}}
\put(64,122){\line(2,1){15.4}}}

\multiput(164,0)(18,-30){4}
{\put(-64,158){\line(-4,-3){15.2}}
\put(-79.2,146.6){\line(0,-1){17}}
\put(-46,128){\line(-3,-1){18}}
\put(-64,122){\line(-2,1){15.4}}}

\put(-54,-90)
{\put(46,128){\line(4,-3){10}}
\put(56,120.5){\line(3,-5){6}}
\put(64,98){\line(-1,6){2.1}}}

\put(218,-90)
{\put(-46,128){\line(-4,-3){10}}
\put(-56,120.5){\line(-3,-5){6}}
\put(-64,98){\line(1,6){2.1}}}

\multiput(-54,-90)(36,0){4}
{\put(64,98){\line(1,6){3.3}}
\put(100,98){\line(-1,6){3.3}}
\put(82,125){\line(2,-1){15}}
\put(82,125){\line(-2,-1){15}}}

\multiput(82,108.5)(18,-30){3}{\circle{30}}
\multiput(64,78.5)(18,-30){2}{\circle{30}}
\multiput(46,48.5)(18,-30){1}{\circle{30}}

\thinlines
\multiput(82,149)(36,0){1}{\line(2,1){18}}
\multiput(82,149)(36,0){1}{\line(-2,1){18}}
\multiput(82,128)(36,0){1}{\line(0,1){21}}
\multiput(46,128)(36,0){2}{\line(2,-1){18}}
\multiput(82,128)(36,0){2}{\line(-2,-1){18}}
\multiput(64,98)(36,0){2}{\line(0,1){21}}
\multiput(28,98)(36,0){3}{\line(2,-1){18}}
\multiput(64,98)(36,0){3}{\line(-2,-1){18}}
\multiput(46,68)(36,0){3}{\line(0,1){21}}
\multiput(10,68)(36,0){4}{\line(2,-1){18}}
\multiput(46,68)(36,0){4}{\line(-2,-1){18}}
\multiput(28,38)(36,0){4}{\line(0,1){21}}
\multiput(-8,38)(36,0){5}{\line(2,-1){18}}
\multiput(28,38)(36,0){5}{\line(-2,-1){18}}
\multiput(10,8)(36,0){5}{\line(0,1){21}}
\put(-26,8){\line(1,0){216}}
\put(-26,8){\line(3,5){108}}
\put(190,8){\line(-3,5){108}}
\end{picture}
\end{center}
\begin{center}
\stepcounter{figure}
Figure \thefigure: $\tilde{\Gamma}$ for the standard simplicial decomposition
\end{center}
which is a fattening of the following graph $\Gamma$. By results in \cite{N}, we can simplectically deform $C$ to symplectic curve $\hat{C}$ such that $F(\hat{C}) = \Gamma$\\
\begin{center}
\setlength{\unitlength}{1.2pt}
\begin{picture}(200,160)(-30,10)
\put(150,48){\line(-3,-5){18}}
\put(6,48){\line(3,-5){18}}
\put(60,138){\line(1,0){36}}
\put(132,78){\line(-3,-5){36}}
\put(24,78){\line(3,-5){36}}
\put(42,108){\line(1,0){72}}
\put(114,108){\line(-3,-5){54}}
\put(42,108){\line(3,-5){54}}
\put(24,78){\line(1,0){108}}
\put(96,138){\line(-3,-5){72}}
\put(60,138){\line(3,-5){72}}
\put(6,48){\line(1,0){144}}
\put(78,168){\line(-3,-5){90}}
\put(78,168){\line(3,-5){90}}
\put(-12,18){\line(1,0){180}}
\multiput(78,168)(36,0){1}{\circle*{4}}
\multiput(60,138)(36,0){2}{\circle*{4}}
\multiput(42,108)(36,0){3}{\circle*{4}}
\multiput(24,78)(36,0){4}{\circle*{4}}
\multiput(6,48)(36,0){5}{\circle*{4}}
\multiput(-12,18)(36,0){6}{\circle*{4}}
\thicklines
\multiput(78,149)(36,0){1}{\line(2,1){18}}
\multiput(78,149)(36,0){1}{\line(-2,1){18}}
\multiput(78,128)(36,0){1}{\line(0,1){21}}
\multiput(42,128)(36,0){2}{\line(2,-1){18}}
\multiput(78,128)(36,0){2}{\line(-2,-1){18}}
\multiput(60,98)(36,0){2}{\line(0,1){21}}
\multiput(24,98)(36,0){3}{\line(2,-1){18}}
\multiput(60,98)(36,0){3}{\line(-2,-1){18}}
\multiput(42,68)(36,0){3}{\line(0,1){21}}
\multiput(6,68)(36,0){4}{\line(2,-1){18}}
\multiput(42,68)(36,0){4}{\line(-2,-1){18}}
\multiput(24,38)(36,0){4}{\line(0,1){21}}
\multiput(-12,38)(36,0){5}{\line(2,-1){18}}
\multiput(24,38)(36,0){5}{\line(-2,-1){18}}
\multiput(6,8)(36,0){5}{\line(0,1){21}}
\end{picture}
\end{center}
\begin{center}
\stepcounter{figure}
Figure \thefigure: $\Gamma$ for the standard simplicial decomposition
\end{center}
If the simplicial decomposition is changed to\\
\begin{center}
\setlength{\unitlength}{1.2pt}
\begin{picture}(200,160)(-30,10)
\put(150,48){\line(-3,-5){18}}
\put(6,48){\line(3,-5){18}}
\put(60,138){\line(1,0){36}}
\put(132,78){\line(-3,-5){36}}
\put(24,78){\line(3,-5){36}}
\put(42,108){\line(1,0){72}}
\put(114,108){\line(-3,-5){54}}
\put(42,108){\line(3,-5){54}}
\put(24,78){\line(1,0){36}}
\put(96,78){\line(1,0){36}}
\put(78,48){\line(0,1){60}}
\put(96,138){\line(-3,-5){72}}
\put(60,138){\line(3,-5){72}}
\put(6,48){\line(1,0){144}}
\put(78,168){\line(-3,-5){90}}
\put(78,168){\line(3,-5){90}}
\put(-12,18){\line(1,0){180}}

\multiput(78,168)(36,0){1}{\circle*{4}}
\multiput(60,138)(36,0){2}{\circle*{4}}
\multiput(42,108)(36,0){3}{\circle*{4}}
\multiput(24,78)(36,0){4}{\circle*{4}}
\multiput(6,48)(36,0){5}{\circle*{4}}
\multiput(-12,18)(36,0){6}{\circle*{4}}
\end{picture}
\end{center}
\begin{center}
\stepcounter{figure}
Figure \thefigure: alternative simplicial decomposition
\end{center}
We have the corresponding $\Gamma$\\
\begin{center}
\setlength{\unitlength}{1.2pt}
\begin{picture}(200,160)(-30,10)
\put(150,48){\line(-3,-5){18}}
\put(6,48){\line(3,-5){18}}
\put(60,138){\line(1,0){36}}
\put(132,78){\line(-3,-5){36}}
\put(24,78){\line(3,-5){36}}
\put(42,108){\line(1,0){72}}
\put(114,108){\line(-3,-5){54}}
\put(42,108){\line(3,-5){54}}
\put(24,78){\line(1,0){36}}
\put(96,78){\line(1,0){36}}
\put(78,48){\line(0,1){60}}
\put(96,138){\line(-3,-5){72}}
\put(60,138){\line(3,-5){72}}
\put(6,48){\line(1,0){144}}
\put(78,168){\line(-3,-5){90}}
\put(78,168){\line(3,-5){90}}
\put(-12,18){\line(1,0){180}}

\multiput(78,168)(36,0){1}{\circle*{4}}
\multiput(60,138)(36,0){2}{\circle*{4}}
\multiput(42,108)(36,0){3}{\circle*{4}}
\multiput(24,78)(36,0){4}{\circle*{4}}
\multiput(6,48)(36,0){5}{\circle*{4}}
\multiput(-12,18)(36,0){6}{\circle*{4}}
\thicklines
\multiput(78,149)(36,0){1}{\line(2,1){18}}
\multiput(78,149)(36,0){1}{\line(-2,1){18}}
\multiput(78,128)(36,0){1}{\line(0,1){21}}
\multiput(42,128)(36,0){2}{\line(2,-1){18}}
\multiput(78,128)(36,0){2}{\line(-2,-1){18}}
\multiput(60,98)(36,0){2}{\line(0,1){21}}
\multiput(24,98)(72,0){2}{\line(2,-1){18}}
\multiput(60,98)(26,-20){2}{\line(1,-2){10}}
\put(70,78){\line(1,0){16}}
\multiput(60,98)(72,0){2}{\line(-2,-1){18}}
\multiput(96,98)(-26,-20){2}{\line(-1,-2){10}}
\multiput(42,68)(72,0){2}{\line(0,1){21}}
\multiput(6,68)(36,0){2}{\line(2,-1){18}}
\multiput(114,68)(36,0){1}{\line(2,-1){18}}
\multiput(42,68)(36,0){1}{\line(-2,-1){18}}
\multiput(114,68)(36,0){2}{\line(-2,-1){18}}
\multiput(24,38)(36,0){4}{\line(0,1){21}}
\multiput(-12,38)(36,0){5}{\line(2,-1){18}}
\multiput(24,38)(36,0){5}{\line(-2,-1){18}}
\multiput(6,8)(36,0){5}{\line(0,1){21}}
\end{picture}
\end{center}
\begin{center}
\stepcounter{figure}
Figure \thefigure: $\Gamma$ for alternative simplicial decomposition
\end{center}
Change the simplicial decomposition further, we will get\\
\begin{center}
\setlength{\unitlength}{1.2pt}
\begin{picture}(200,160)(-30,10)
\put(61,136){\line(2,-1){52}}
\put(150,48){\line(-3,-5){18}}
\put(6,48){\line(3,-5){18}}
\put(60,138){\line(1,0){36}}
\put(132,78){\line(-3,-5){36}}
\put(24,78){\line(3,-5){36}}
\put(42,108){\line(1,0){72}}
\put(114,108){\line(-3,-5){54}}
\put(42,108){\line(3,-5){54}}
\put(24,78){\line(1,0){36}}
\put(96,78){\line(1,0){36}}
\put(78,48){\line(0,1){60}}
\put(78,108){\line(-3,-5){54}}
\put(60,138){\line(3,-5){72}}
\put(6,48){\line(1,0){144}}
\put(78,168){\line(-3,-5){90}}
\put(78,168){\line(3,-5){90}}
\put(-12,18){\line(1,0){180}}
\multiput(78,168)(36,0){1}{\circle*{4}}
\multiput(60,138)(36,0){2}{\circle*{4}}
\multiput(42,108)(36,0){3}{\circle*{4}}
\multiput(24,78)(36,0){4}{\circle*{4}}
\multiput(6,48)(36,0){5}{\circle*{4}}
\multiput(-12,18)(36,0){6}{\circle*{4}}
\end{picture}
\end{center}
\begin{center}
\stepcounter{figure}
Figure \thefigure: another alternative simplicial decomposition
\end{center}
We have the corresponding $\Gamma$\\
\begin{center}
\setlength{\unitlength}{1.2pt}
\begin{picture}(200,160)(-30,10)
\put(61,136){\line(2,-1){52}}
\put(150,48){\line(-3,-5){18}}
\put(6,48){\line(3,-5){18}}
\put(60,138){\line(1,0){36}}
\put(132,78){\line(-3,-5){36}}
\put(24,78){\line(3,-5){36}}
\put(42,108){\line(1,0){72}}
\put(114,108){\line(-3,-5){54}}
\put(42,108){\line(3,-5){54}}
\put(24,78){\line(1,0){36}}
\put(96,78){\line(1,0){36}}
\put(78,48){\line(0,1){60}}
\put(78,108){\line(-3,-5){54}}
\put(60,138){\line(3,-5){72}}
\put(6,48){\line(1,0){144}}
\put(78,168){\line(-3,-5){90}}
\put(78,168){\line(3,-5){90}}
\put(-12,18){\line(1,0){180}}
\multiput(78,168)(36,0){1}{\circle*{4}}
\multiput(60,138)(36,0){2}{\circle*{4}}
\multiput(42,108)(36,0){3}{\circle*{4}}
\multiput(24,78)(36,0){4}{\circle*{4}}
\multiput(6,48)(36,0){5}{\circle*{4}}
\multiput(-12,18)(36,0){6}{\circle*{4}}
\thicklines
\multiput(78,149)(36,0){1}{\line(2,1){18}}
\multiput(78,149)(36,0){1}{\line(-2,1){18}}
\multiput(78,149)(4,-30){2}{\line(2,-3){14}}
\multiput(114,129)(-32,-10){2}{\line(-1,0){23}}
\multiput(82,119)(-27,-10){1}{\line(1,1){10}}
\multiput(42,128)(36,0){1}{\line(2,-1){18}}
\multiput(78,128)(36,0){0}{\line(-2,-1){18}}
\multiput(60,98)(36,0){1}{\line(0,1){21}}
\multiput(24,98)(72,0){2}{\line(2,-1){18}}
\multiput(60,98)(26,-20){2}{\line(1,-2){10}}
\put(70,78){\line(1,0){16}}
\multiput(60,98)(72,0){2}{\line(-2,-1){18}}
\multiput(96,98)(-26,-20){2}{\line(-1,-2){10}}
\multiput(42,68)(72,0){2}{\line(0,1){21}}
\multiput(6,68)(36,0){2}{\line(2,-1){18}}
\multiput(114,68)(36,0){1}{\line(2,-1){18}}
\multiput(42,68)(36,0){1}{\line(-2,-1){18}}
\multiput(114,68)(36,0){2}{\line(-2,-1){18}}
\multiput(24,38)(36,0){4}{\line(0,1){21}}
\multiput(-12,38)(36,0){5}{\line(2,-1){18}}
\multiput(24,38)(36,0){5}{\line(-2,-1){18}}
\multiput(6,8)(36,0){5}{\line(0,1){21}}
\end{picture}
\end{center}
\begin{center}
\stepcounter{figure}
Figure \thefigure: the corresponding $\Gamma$
\end{center}

\subsection{Lagrangian fibration with codimension 2 singular locus}
With all this preparation, now we would like to address the meaning of {\bf near the large complex limit}. Consider the Calabi-Yau hypersurface $X_s$ in $P_\Sigma$ defined by\\
\[
s = \sum_{m(\not=0)\in\Delta}a_ms_m + \psi s_0=0
\]\\
Let $\Delta^0$ denote the set of integral points in the 2-skeleton of $\Delta$. Or in another word, the integral points in $\Delta$ that is not in the interior of 3-faces and also not the origion. For $m\in \Delta$, let $|a_m| = e^{-w_m}$, $w'_m = w_m-w_{0}$ (recall $\psi = a_{0}$). We have two functions $w=(w_m)_{m\in \Delta}$, $w'=(w'_m)_{m\in \Delta}$ defined on the integral points in $\Delta$.\\
\begin{de}
\label{fd}
$w'=(w'_m)_{m\in \Delta}$ is called convex with respect to  $\Delta^0$, if for any $m' \in \Delta^0$ there exist a linear function $n$ such that $n(m')= w'_{m'}$ and $n(m)\leq w_m$ for $m\in \Delta^0\backslash \{m'\}$.\\
\end{de}
The two convexities we defined are closely related. The 2-skeleton\\
\[
\Delta^0 = \bigcup_{\sigma\in \Sigma(2)} \Delta_{\sigma^*}.
\]\\
Let\\
\[
w_\sigma = w|_{\Delta_{\sigma^*}}.
\]\\
Then we have the following lemma\\
\begin{lm}
$w'=(w'_m)_{m\in \Delta}$ is convex with respect to  $\Delta^0$ for $-w_{m_0}$ large enough, if and only if $w_\sigma$ is convex on $\Delta_{\sigma^*}$ in the sense of definition \ref{fb} for all $\sigma\in \Sigma(2)$.\\
\end{lm}
With this lemma in mind, we can make the concept of near the large complex limit more precise as follows\\
\begin{de}
\label{fe}
The quintic Calabi-Yau hypersurface $X_s$ is said to be {\bf near the large complex limit}, if $w'=(w'_m)_{m\in \Delta}$ is convex with respect to  $\Delta^0$ and $|w|$ is large.\\
\end{de}
When $X_s$ is near the large complex limit, the corresponding $w_\sigma$ is convex on $\Delta_{\sigma^*}$. When $w$ is generic, by previous construction, $w_\sigma$ determine a 1-dimensional graph $\Gamma_\sigma$ on $\Delta_{\sigma^*}$. Let\\
\[
\Gamma = \bigcup_{\sigma\in \Sigma(2)} (\hat{\pi}_v^{(2)})^{-1}(\Gamma_\sigma)
\]\\
By theorem \ref{fc}, $\tilde{\Gamma}_\sigma = F(C_\sigma)$ is a fattening of $\Gamma_\sigma$. Therefore, 
\[
\tilde{\Gamma} = \bigcup_{\sigma\in \Sigma(2)} (\hat{\pi}_v^{(2)})^{-1}(\tilde{\Gamma}_\sigma)
\]\\
is a fattening of 1-dimensional graph $\Gamma$. Let $\Gamma = \Gamma^1\cup \Gamma^2\cup \Gamma^3$, where $\Gamma^1$ is the smooth part of $\Gamma$, $\Gamma^2$ is the singular part of $\Gamma$ that will map to the interior of the 2-skeleton of $\Delta$ under $\hat{\pi}_v$, $\Gamma^3$ is the singular part of $\Gamma$ that will map to the 1-skeleton of $\Delta$ under $\hat{\pi}_v$. Apply second part of the theorem \ref{fc}, also with help of some other result in \cite{N}, we can produce a Lagrangian fibration $\hat{F}: X_\infty \rightarrow \partial \Delta_v$ such that $\hat{F}(C)=\Gamma$. Then we have\\
\begin{th}
\label{ff}
When $X_\psi$ is near the large complex limit, start with Lagrangian fibration $\hat{F}$ the flow of $V$ will produce a Lagrangian fibration $\hat{F}_s: X_s \rightarrow \partial\Delta_v$. There are 4 types of fibres.\\
(i). For $p\in \partial\Delta_v\backslash \Gamma$, $\hat{F}_s^{-1}(p)$ is a smooth Lagrangian 3-torus.\\
(ii). For $p\in \Gamma^1$, $\hat{F}_s^{-1}(p)$ is a type I singular fibre.\\ 
(iii). For $p\in \Gamma^2$, $\hat{F}_s^{-1}(p)$ is a type II singular fibre.\\
(iv). For $p\in \Gamma^3$, $\hat{F}_s^{-1}(p)$ is a type III singular fibre.\\
\end{th}
{\bf Remark:} Although the Lagrangian fibration produced by this theorem have exactly the same topological structure as the expected special Lagrangian fibration, this fibration map is still not $C^{\infty}$ map (merely Lipschitz). But I believe that a small pertubation of this map will give a $C^{\infty}$ Lagrangian fibration with the same topological structure. In \cite{lag2}, we were able to make the map $C^\infty$ away from $C$. We also indicated how to modify singular fibres around $C$. Yet, we still fall short to completely make the fibration $C^\infty$. In any case, we do not really use $C^\infty$ property of the fibration map in the following discussion.\\\\

\se{Duality map of convex polyhedrons}
As we know, $(p_v,p_w)$ parametrizes \k moduli and complex moduli of Calabi-Yau hypersurface $X$ in $P_{\Sigma^v}$. From symmetric point of view, $(p_v,p_w)$ also parametrizes complex moduli and \k moduli of the mirror Calabi-Yau hypersurface $Y$ in $P_{\Sigma^w}$. According to the previous section, we can construct Lagrangian torus fibrations\\
\[
\hat{F}_{w,v}: X \rightarrow \partial \Delta_v,\ \ \ \hat{F}_{v,w}: Y \rightarrow \partial \Delta_w.
\]\\
with singular locuses $\Gamma_{w,v}$ and $\Gamma_{v,w}$.\\\\
One key ingredient of SYZ mirror conjecture is an identification of the bases of the fibrations $\phi: \Delta_w \rightarrow \Delta_v$, such that $\phi(\Gamma_{v,w})=\Gamma_{w,v}$. It turns out that this identification can be constructed purely combinatorically depending on the two convex piecewise linear convex functions $p_v,p_w$, without refering to symplectic form, Calabi-Yau etc. We will descuss this purely combinatoric construction in this section.\\\\
Assume that $\Delta$ is a reflexive polyhedron, $\Sigma$ the corresponding fan. Faces of $\Delta$ ($\Delta^\vee$) is in 1-1 correspondence with cones in $\Sigma^\vee$ ($\Sigma$). For a cone $\sigma$ in $\Sigma^\vee$ ($\Sigma$), let $\Delta_\sigma = \sigma\cap (\Delta\backslash\{0\})$ ($(\Delta^\vee)_{\sigma} = \sigma\cap (\Delta^\vee\backslash\{0\})$) be the corresponding face of $\Delta$ ($\Delta^\vee$). This correspondence preserve the inclusion relation. There is a natural dual correspondence between faces of $\Delta$ and $\Delta^\vee$. For $\Delta_\sigma$ a face of $\Delta$, the dual face\\
\[
(\Delta_\sigma)^* = \{n\in \Delta^\vee|\langle m,n\rangle =-1,\ {\rm for}\ m\in \Delta_\sigma\}.
\]\\
is a face of $\Delta^\vee$. Then there exist cone $\sigma^*\in \Sigma$ such that $(\Delta_\sigma)^* = \Delta^\vee_{\sigma^*}$ In this way, there is a natural identification of faces of $\Delta$ and $\Delta^\vee$, or equivalently, a natural identification of the sets $\Sigma^\vee$ and $\Sigma$. These identifications reverse the inclusion relation, which is a partial relation.\\\\
A simplex in the baricenter subdivision of $\Delta$ ($\Delta^\vee$) can be expressed by an inclusion related sequence of cones $\{\sigma_k\}_{k=0}^l$ in $\Sigma^\vee$ ($\Sigma$) that satisfies $\sigma_k\subset \sigma_{k+1}$. Here $\sigma_k$ is best understood to represent the baricenter of $\Delta_{\sigma_k}$. $\hat{\sigma} = \{\sigma_k\}_{k=0}^l$ denotes the $l$-simplex $\hat{\sigma}$ with vertices being the baricenters of $\Delta_{\sigma_k}$ for $k=0,\cdots,l$.\\
\begin{prop}
\label{ga}
The dual correspondence induces naturally a piecewise linear homeomorphism $\Delta \rightarrow \Delta^\vee$ with respect to the baricenter subdivision. Where baricenter of $\Delta_\sigma$ is map to baricenter of $\Delta^\vee_{\sigma^*}$, and simplex $\hat{\sigma} = \{\sigma_k\}_{k=0}^l$ is mapped to $\hat{\sigma}^* = \{\sigma^*_k\}_{k=0}^l$ linearly.\\
\end{prop}
From now on, we will introduce the notation $\hat{\Sigma}^\vee$ ($\hat{\Sigma}$) (not as a fan) to index simpleces $\hat{\sigma} = \{\sigma_k\}_{k=0}^l$ ($\hat{\sigma}^* = \{\sigma^*_k\}_{k=0}^l$) in the baricenter subdivision of $\Delta$ ($\Delta^\vee$).\\\\   
Assume that $p_w$ ($p_v$) is a convex piecewise linear function on $M$ ($N$) with respect to $\Delta$ ($\Delta^\vee$). $p_w(m)=-w_m$ ($p_v(n)=-v_n$) for $m\in\Delta\backslash\{0\}$ ($n\in\Delta^\vee\backslash\{0\}$). $p_w$ ($p_v$) naturally determines a fan $\Sigma^w$ ($\Sigma^v$) for $M$ ($N$) that is compatible with $\Delta^\vee$ ($\Delta$). $p_w$ ($p_v$) also naturally determines a polyhedron $\Delta^w\subset N$ ($\Delta^v\subset M$) consists of $n\in N$ ($m\in M$) that as a linear function on $M$ ($N$) is greater or equal to $p_w$ ($p_v$). In particular, the fan $\Sigma$ ($\Sigma^\vee$) of the anti-canonical model $P_{\Sigma}$ ($P_{\Sigma^\vee}$) is determined by $p_\Delta$ ($p_{\Delta^\vee}$), $p_\Delta(n)=-1$ ($p_{\Delta^\vee}(m)=-1$) for $n\in{\Delta^\vee} \backslash\{0\}$ ($m\in\Delta\backslash\{0\}$), for which the corresponding polyhedron is $\Delta$ ($\Delta^\vee$).\\\\
$p_w$ ($p_v$) determines a polyhedron decomposition $Z^w$ ($Z^v$) of $\Delta$ ($(\Delta^\vee)$). Let $\hat{Z}^w$ ($\hat{Z}^v$) be the baricenter subdivision of $Z^w$ ($Z^v$).\\\\
Restrict to each face, for $\sigma\in\Sigma^\vee$, $p_w$ ($p_v$) determines a polyhedron decomposition $Z^w_\sigma$ ($Z^v_{\sigma^*}$) of $\Delta_\sigma$ ($(\Delta^\vee)_{\sigma^*}$). Here $Z^w_\sigma$ ($Z^v_{\sigma^*}$) only consists of those open polyhedrons whose baricenters are in the interior of $\Delta_\sigma$ ($(\Delta^\vee)_{\sigma^*}$). We have\\
\[
Z^w = \bigcup_{\sigma\in \Sigma} Z^w_{\sigma^*} 
\]\\
\[
Z^v = \bigcup_{\sigma\in \Sigma} Z^v_{\sigma} 
\]\\
Let $\hat{Z}^w_\sigma$ ($\hat{Z}^v_{\sigma^*}$) be the baricenter subdivision of $Z^w_\sigma$ ($Z^v_{\sigma^*}$), and simplices in $\hat{Z}^w_\sigma$ ($\hat{Z}^v_{\sigma^*}$) have vertices in $Z^w_\sigma$ ($Z^v_{\sigma^*}$).\\\\
For $\hat{\sigma} = \{\sigma_k\}_{k=0}^l \in\hat{\Sigma}^\vee$ ($\hat{\sigma}^* = \{\sigma^*_k\}_{k=0}^l \in \hat{\Sigma}$), define\\
\[
\hat{Z}^w_{\hat{\sigma}} = \{ \hat{\alpha} =\{\hat{\alpha}_0,\cdots, \hat{\alpha}_l\}|\hat{\alpha}_k\in \hat{Z}^w_{\sigma_k},\   \hat{\alpha}_k < \hat{\alpha}_{k+1},\ {\rm for\ all}\ k.\}
\]\\
\[
\hat{Z}^v_{\hat{\sigma}^*} = \{ \hat{\alpha} =\{\hat{\alpha}_0,\cdots, \hat{\alpha}_l\}|\hat{\alpha}_k\in \hat{Z}^v_{\sigma^*_k},\   \hat{\alpha}_k < \hat{\alpha}_{k-1},\ {\rm for\ all}\ k.\}
\]\\
where $\hat{\alpha} < \hat{\beta}$ means that each polyhedron in $\hat{\alpha}$ is a face of every polyhedron in $\hat{\beta}$. Then\\
\[
\hat{Z}^w = \bigcup_{\hat{\sigma}\in \hat{\Sigma}} \hat{Z}^w_{\hat{\sigma}^*} 
\]\\
\[
\hat{Z}^v = \bigcup_{\hat{\sigma}\in \hat{\Sigma}} \hat{Z}^v_{\hat{\sigma}} 
\]\\
Now we are ready to discuss the polyhedron decomposition on $\Delta^v$. Consider the piecewise linear map $\pi_v: \Delta^v\rightarrow \Delta$, for $\sigma\in \Sigma^\vee$, $\pi_v$ over the interior of $\Delta_\sigma$ is a fibration with compact fibres. A simplex in a fibre is corresponding to a simplex of complimentary dimension in $Z^v_{\sigma^*}$. In this way, we may identify the baricenter subdivision of a fibre by $\hat{Z}^v_{\sigma^*}$. (the reason for baricenter subdivision is that after it the correspondence will preserve the dimension of the simplex, instead of go to the complimentary dimension.) On the other hand, $\hat{Z}^w_\sigma$ form a simplicial decomposition of interior of $\Delta_\sigma$. For $\hat{\alpha} \in \hat{Z}^w_\sigma$, $\hat{\beta} \in \hat{Z}^v_{\sigma^*}$, $\pi_v^{-1}(\hat{\alpha})\cap\hat{\beta} \cong \hat{\alpha} \times \hat{\beta}$ form a polyhedron $\Delta^v_{\hat{\alpha},\hat{\beta}}$ in $\Delta^v$. Let $\hat{Z}^{v,w}_\sigma \cong \hat{Z}^w_\sigma \times \hat{Z}^v_{\sigma^*}$ represent the set of this kind of polyhedron in $\Delta^v$. Then for $\hat{\sigma} = \{\sigma_k\}_{k=0}^l \in\hat{\Sigma}^\vee$, define\\
\[
\hat{Z}^{v,w}_{\hat{\sigma}} = \{ \hat{\gamma} =\{\hat{\alpha}_0\times \hat{\beta}_0,\cdots, \hat{\alpha}_l\times \hat{\beta}_l\}|\hat{\alpha}_k\times \hat{\beta}_k\in \hat{Z}^{v,w}_{\sigma_k},\   \hat{\alpha}_k < \hat{\alpha}_{k+1},\ \hat{\beta}_{k+1} < \hat{\beta}_k\ {\rm for\ all}\ k.\}
\]\\
\[
\hat{Z}^{v,w} = \bigcup_{\hat{\sigma}\in \hat{\Sigma}^\vee} \hat{Z}^{v,w}_{\hat{\sigma}} 
\]\\
Then $\hat{Z}^{v,w}$ forms a polyhedron decomposition of $\Delta^v$.\\
\[
\hat{\gamma}' =\{\hat{\alpha}_0'\times \hat{\beta}_0',\cdots, \hat{\alpha}_{l'}'\times \hat{\beta}_{l'}'\} \in \hat{Z}^{v,w}_{\hat{\sigma}'}
\]\\
is a face of\\
\[
\hat{\gamma} =\{\hat{\alpha}_0\times \hat{\beta}_0,\cdots, \hat{\alpha}_l\times \hat{\beta}_l\}\in \hat{Z}^{v,w}_{\hat{\sigma}}
\]\\
if and only if ${\hat{\sigma}'}\subset {\hat{\sigma}}$ (more precisely, $l'\leq l$, $\sigma_k'=\sigma_k$ for $1\leq k\leq l'$) and ${\hat{\alpha}_k'}\subset {\hat{\alpha}_k}$, ${\hat{\beta}_k'}\subset {\hat{\beta}_k}$ for $1\leq k\leq l'$.\\\\
In completely symmetric way, we can define\\
\[
\hat{Z}^{w,v} = \bigcup_{\hat{\sigma}\in \hat{\Sigma}} \hat{Z}^{w,v}_{\hat{\sigma}} 
\]\\
as a polyhedron decomposition of $\Delta^w$. Since\\
\[
\hat{Z}^{v,w}_\sigma \cong \hat{Z}^w_\sigma \times \hat{Z}^v_{\sigma^*} \cong \hat{Z}^{w,v}_{\sigma^*}
\]\\
and $\sigma \rightarrow \sigma^*$ naturally identify $\Sigma$ with $\Sigma^\vee$, we have\\
\begin{prop}
\[
\hat{Z}^{w,v} \cong \hat{Z}^{v,w}
\]
and subface is identified with subface.\\
\end{prop}
This identification naturally induce\\
\begin{th}
\label{gb}
There is a natural piecewise linear homeomorphism\\
\[
\Delta^w \rightarrow \Delta^v
\]
with respect to polyhedron decompositions $\hat{Z}^{w,v} \cong \hat{Z}^{v,w}$.\\
\end{th}
There is another way to see the theorem. for $\sigma\in \Sigma^\vee$, $\pi_v$ over the interior of $\Delta_\sigma$ is a fibration with compact fibres. A simplex $\beta^*$ in a fibre is corresponding to a simplex of complimentary dimension $\beta$ in $Z^v_{\sigma^*}$. On the other hand, $Z^w_\sigma$ form a symplicial decomposition of interior of $\Delta_\sigma$. For $\alpha \in Z^w_\sigma$, $\beta \in Z^v_{\sigma^*}$, the closure of $\pi_v^{-1}(\alpha)\cap\beta^* \cong \alpha \times \beta^*$ in $\Delta^v$ forms a polyhedron $\Delta^v_{\alpha,\beta}$ in $\Delta^v$. Let $Z^{v,w}_\sigma \cong Z^w_\sigma \times Z^v_{\sigma^*}$ represent the set of this kind of polyhedron in $\Delta^v$. Then\\
\[
Z^{v,w} = \bigcup_{\sigma\in \Sigma^\vee} Z^{v,w}_{\sigma} 
\]\\
forms a polyhedron decomposition of $\Delta^v$.\\
\[
\alpha' \times (\beta')^* \in Z^{v,w}_{\sigma'} \cong Z^w_{\sigma'} \times Z^v_{(\sigma')^*} 
\]\\
is a subface of\\
\[
\alpha \times \beta \in Z^{v,w}_\sigma \cong Z^w_\sigma \times Z^v_{\sigma^*} 
\]\\
if and only if $\sigma'$, $\alpha'$, $\beta'$ is a subface of $\sigma$, $\alpha$, $\beta$.\\\\
In completely symmetric way, we can define $Z^{w,v}_\sigma \cong Z^v_\sigma \times Z^w_{\sigma^*}$ and\\
\[
Z^{w,v} = \bigcup_{\sigma\in \Sigma^\vee} Z^{w,v}_{\sigma} 
\]\\
forms a polyhedron decomposition of $\Delta^w$.\\
The natural identification of $Z^{v,w}$ and $Z^{w,v}$ is a duality relation (inclusion relation is reversed). Then in the spirit of proposition \ref{ga}, this dual relation induce a idnetification $\Delta^w \rightarrow \Delta^v$ that is piecewise linear with respect to the baricenter subdivisons $\hat{Z}^{v,w}$ and $\hat{Z}^{w,v}$.\\\\

\se{SYZ mirror construction for Calabi-Yau hypersurface}
To fully establish SYZ construction, we also need to establish the duality relation of the Lagrangian torus fibres in the Lagrangian fibrations of Calabi-Yau hypersurface and its mirror. There are many way to see this, for example, one may compute monodromy operators of the two fibrations, and show that they are dual to each other. We will use a more direct method, we will establish a canonical characterization of Lagrangian torus fibre of the Lagrangian fibrations of Calabi-Yau hypersurface in $P_{\Sigma_v}$ and its mirror in $P_{\Sigma_w}$, and prove that under this canonical characterization the two fibres related via the base identification $\phi$ constructed in the previous section are canonically dual to each other.\\\\
Recall the \k moduli (the movable cone) of Calabi-Yau hypersurface $Y\subset P_{\Sigma^w}$ can be characterized as\\
\[
\tau = \bigcup_{Z\in \tilde{Z}} \tau_Z
\]\\
where $\tilde{Z}$ denotes the set of simplicial decomposition of $\Delta^0$. For a particular simplicial decomposition $Z\in \tilde{Z}$, $\tau_Z$ is the set of $w=(w_m)_{m\in \Delta^0}$ (understood as a piecewise linear function on $M$) that is convex with respect to $\Delta^0$ in the sense of definition \ref{fd}, and determine the simplicial decomposition $Z$ of $\Delta^0$, modulo linear function on $M$.\\\\
The monomial divisor map gives us a natural identification of the complex moduli of Calabi-Yau hypersurface $X\subset P_{\Sigma^v}$ near the large complex limit with the complexified \k moduli of Calabi-Yau hypersurface $Y\subset P_{\Sigma^w}$, $\tau_{\bf C} = (\tilde{N}\otimes_{\bf Z}{\bf R} + i\tau)/\tilde{N}$.\\\\
For any $u=(u_m)_{m\in \Delta^0}\in (\tilde{N}\otimes_{\bf Z}{\bf R} + i\tau)/\tilde{N}$, consider the Calabi-Yau hypersurface $X_u\subset P_{\Sigma^v}$ defined by\\
\[
s_u(z)= \sum_{m\in \Delta^0} e^{2\pi iu_m}s_m + s_0=0.
\]\\
Let $w_m = {\rm Im}(u_m)$, then $w=(w_m)_{m\in \Delta^0}\in \tau$.Lagrangian torus fibration $X_u\rightarrow \partial \Delta_v$ is constructed by deforming the natural Lagrangian torus fibration of the large complex limit $X_\infty$ via gradient flow. Where $X_\infty$ is defined by $s_0=0$.\\\\
For any integral $n \in \partial\Delta^\vee\subset N$, there is a corresponding dimension 3 face $\alpha_n$ of $\Delta_v$ defined as\\
\[
\alpha_n = \{m\in \Delta_v|v_n\langle m,n\rangle =-1.\}
\]\\
Namely $v_nn$ is the unique supporting function of $\alpha_n$. Clearly, fibres of the fibration $X_\infty\rightarrow \partial \Delta_v$ over $\alpha_n^0$ (interior of $\alpha_n$) are naturally identified with\\
\[
T_n \cong (N_n \otimes_{\bf Z} {\bf R})/N_n
\]\\
where $N_n = N/\{{\bf Z}n\}$. Since Lagrangian torus fibration $X_u\rightarrow \partial \Delta_v$ is a deformation of fibration $X_\infty\rightarrow \partial \Delta_v$, we have\\
\begin{prop}
3-torus fibres of the Lagrangian fibration $X_u\rightarrow \partial \Delta_v$ over interior of $\alpha_n$ can be naturally identified with $T_n$.\\
\end{prop}
The dual torus of $T_n$ is naturally\\
\[
T_n^\vee \cong (M_n \otimes_{\bf Z} {\bf R})/M_n
\]\\
where $M_n = n^\perp \subset M$.\\\\
Similarly, for any integral $m\in \partial\Delta$, we can define\\
\[
T_m^\vee \cong (N_m \otimes_{\bf Z} {\bf R})/N_m
\]\\
where $N_m = m^\perp \subset N$. We have\\
\begin{prop}
\label{ha}
3-torus fibre $X_b$ of the Lagrangian fibration $X_u\rightarrow \partial \Delta_v$ over $b$ in a small neighborhood $U_m$ of $\hat{\pi}_v^{-1}(m)\in \partial \Delta_v$ can be naturally identified with $T_m^\vee$. In addition, if $b\in \alpha_n^0$, then the following diagram commutes\\
\[
\begin{array}{ccc}
X_b&\rightarrow& T_m^\vee\\
\downarrow&&\downarrow\\
T_n&\leftarrow& T_{\bf R}\\
\end{array}
\]\\
\end{prop}
{\bf Proof:}
$g\in T$ naturally acts on $P_{\Sigma_v}$, and also on the space of Calabi-Yau hypersurfaces in $P_{\Sigma_v}$. The action acturally preserve our slice. Under the action of $g$, Calabi-Yau hypersurfaces $X_u$ defined by $s_u(z)=0$ is mapped to $X_{g(u)}$ defined by $s_{g(u)}(z)=s_u(g^{-1}(z))=0$. Recall that the moment map\\
\[
F_{w,v}(x) = \int_{m\in\Delta_v} \frac{|s_m|_w^2}{|s|_{w,v}^2}m,
\]\\
maps $P_{\Sigma_v}$ to $\Delta_v$. The action of $g$ on $u$ will naturally induce an action on $w$. Clearly, the image $F_{w,v}(X_u) \subset \Delta_v$ is invariant under the action of $g$. More precisely, the following diagram commutes.\\
\begin{center}
\setlength{\unitlength}{1.3pt}
\begin{picture}(100,50)(0,0)
\put(30,40){\makebox(0,0){$X_u$}}
\put(80,40){\makebox(0,0){$X_{g(u)}$}}
\put(30,10){\makebox(0,0){$\Delta_v$}}
\put(80,10){\makebox(0,0){$\Delta_v$}}
\put(42,40){\vector(1,0){26}}
\put(28,32){\vector(0,-1){15}}
\put(42,10){\vector(1,0){26}}
\put(78,32){\vector(0,-1){15}}
\put(55,45){\makebox(0,0){\footnotesize{$g$}}}
\put(55,15){\makebox(0,0){\footnotesize{$ $}}}
\put(37,25){\makebox(0,0){\footnotesize{$F_{w,v}$}}}
\put(91,25){\makebox(0,0){\footnotesize{$F_{g(w),v}$}}}
\end{picture}
\end{center}
Think of $w$ as a convex piecewise linear function on $M$, the action of $g\in T$ is just modifying $w$ by linear function. Since $w$ is strongly convex, for integral $m'\in \partial \Delta$, by suitable modification of linear function, we may assume that $u_{m'}=0$ and $w_m>0$ for $m\in \partial \Delta\backslash \{m'\}$. Then\\
\[
s_u(z)= s_0 + s_{m'} + \sum_{m\in \Delta_0\backslash \{m'\}} e^{2\pi iu_m}s_m.
\]\\
Let $s_u =s_0 + \tilde{s}_u$, then\\
\[
\tilde{s}_u(z)= s_{m'} + \sum_{m\in \Delta_0\backslash \{m'\}} e^{2\pi iu_m}s_m. 
\]\\
For the purpose of gradient flow, consider meromorphic function\\
\[
q= \frac{s_0(z)}{\tilde{s}_u(z)}
\]\\
$\nabla f$ for $f={\rm Re}(q)$ is our gradient vector field. We actually use the flow of $V = \frac{\nabla f}{|\nabla f|^2}$.\\\\ 
When $s_u(z)$ is near the large complex limit, for $z$ near $F_{w,v}^{-1}\circ \hat{\pi}_v^{-1}(m')$, we have $\tilde{s}_u(z)= s_{m'} + o(1)$. Hence, $\tilde{s}_u(z)$ is non-vanishing near $F_{w,v}^{-1}\circ \hat{\pi}_v^{-1}(m')$, and the flow of $V$ is moving Lagrangian 3-torus fibre near $F_{w,v}^{-1}\circ \hat{\pi}_v^{-1}(m')$ entirely away from $X_\infty$ into the large torus $T$. We need to show that this 3-torus in $T$ is of the same homotopy class as $T_{m'}^\vee \subset T$.\\\\
In general, a face $\alpha_n$ of $\Delta_v$ intersects $U_{m'}$ is and only if $\langle m',n\rangle=-1$. Under this condition, the composition of\\
\begin{equation}
\label{hb}
T_{m'}^\vee \rightarrow T_{\bf R}\rightarrow T_n 
\end{equation}\\
is an isomorphism. This condition also ensure that we can choose $T$-invariant coordinate $x=(x_1,x_2,x_3,x_4)$ near $F_{w,v}^{-1}(\alpha_n)$ such that $F_{w,v}^{-1}(\alpha_n)$ is defined by $x_4=0$ and\\
\[
q= \frac{s_0(z)}{\tilde{s}_u(z)} = \frac{s_0(z)}{s_{m'}(z)}(1+o(1)) = x_4(1+o(1)).
\]\\
Namely, $x_4$ is the $T$-invariant function on $T$ corresponding to $-m'\in M$. Let $x_k=r_ke^{i\theta_k}$ for $i=1,\cdots,4$, then under this coordinate (\ref{hb}) can be expressed as\\
\[
T_m^\vee=\{(\theta_1, \theta_2, \theta_3, 0)\} \rightarrow T_{\bf R}=\{(\theta_1, \theta_2, \theta_3, \theta_4)\}\rightarrow T_n=\{(\theta_1, \theta_2, \theta_3)\} 
\]\\
Start with Lagrangian 3-torus\\
\[
L_0 = \{x||x_i|=r_i, {\rm for}\ i=1,2,3,\ x_4=0\}
\]\\
The flow of $\nabla f$ in $\epsilon$ time will deform $L_0$ to approximately\\
\[
L_\epsilon = \{x||x_i|=r_i, {\rm for}\ i=1,2,3,\ x_4=\epsilon \}
\]\\
This 3-torus is clearly naturally identified with $T_{m'}^\vee$. With this explicit description of $X_b$, it is easy to see that the diagram in the proposition commute.
\begin{flushright} $\Box$ \end{flushright}
Now, let's consider the mirror Calabi-Yau $Y\subset P_{\Sigma^w}$. Also by gradient flow method, we can construct Lagrangian fibration $Y \rightarrow \partial \Delta_w$. For any $m \in \partial \Delta \cap M$, there is a corresponding dimension 3 face $\alpha_m$ of $\Delta_w$ defined as\\
\[
\alpha_m = \{n\in \Delta_w|w_m\langle m,n\rangle =-1.\}
\]\\
Namely $w_mm$ is the unique supporting function of $\alpha_m$. Clearly, fibres of the fibration $Y_\infty\rightarrow \partial \Delta_w$ over $\alpha_m^0$ (interior of $\alpha_m$) are naturally identified with\\
\[
T_m \cong (M_m \otimes_{\bf Z} {\bf R})/M_m
\]\\
where $M_m = M/\{{\bf Z}m\}$. Since Lagrangian torus fibration $Y\rightarrow \partial \Delta_w$ is a deformation of fibration $Y_\infty\rightarrow \partial \Delta_w$, we have\\
\begin{prop}
3-torus fibres of the Lagrangian fibration $Y \rightarrow \partial \Delta_w$ over interior of $\alpha_m$ can be naturally identified with $T_m$.\\
\end{prop}
Recall the natural map $\hat{\pi}_w: \partial \Delta_w \rightarrow \partial \Delta^\vee$. For any integral point $n$ in $\partial\Delta^\vee$, we have\\
\begin{prop}
\label{hc}
3-torus fibre $Y_b$ of the Lagrangian fibration $Y \rightarrow \partial \Delta_w$ over $b$ in a small neighborhood $U_n$ of $\hat{\pi}_w^{-1}(n) \subset \partial \Delta_w$ can be naturally identified with $T_n^\vee$. In addition, if $b\in \alpha_m^0$, then the following diagram commutes\\
\[
\begin{array}{ccc}
Y_b&\rightarrow& T_n^\vee\\
\downarrow&&\downarrow\\
T_m&\leftarrow& T_{\bf R}\\
\end{array}
\]\\
\end{prop}
{\bf Proof:} This is just the mirror statement of proposition \ref{ha}.
\begin{flushright} $\Box$ \end{flushright}
Recall from the theorem \ref{gb}, we have the natural piecewise linear identification of the two base of the Lagrangian fibrations\\
\[
\phi: \partial \Delta_w \rightarrow \partial \Delta_v 
\]\\
Let $\alpha_n^0$, $\alpha_m^0$ denote the interior of $\alpha_n$, $\alpha_m$. Then we have\\
\begin{prop}
$U_m$, $U_n$ can be suitablly chosen such that\\
\[
\phi(U_n) = \alpha_n^0,\ \ \phi(\alpha_m^0) = U_m.
\]
And\\
\[
\partial \Delta_w \backslash \Gamma' = \left(\bigcup_{m\in \Delta^0} \alpha_m^0\right) \cup \left(\bigcup_{n\in \partial \Delta^\vee} U_n\right)
\]\\
\[
\partial \Delta_v \backslash \Gamma = \left(\bigcup_{n\in \partial \Delta^\vee} \alpha_n^0\right) \cup \left(\bigcup_{m\in \Delta^0} U_m\right)
\]\\
\end{prop}
We are now ready to establish the dual relation of the fibres. For any $b\in \partial \Delta_w \backslash \Gamma'$, let $Y_b$ ($X_{\phi(b)}$) denote the fibre of the Lagrangian fibration $Y \rightarrow \partial \Delta_w$ ($X \rightarrow \partial \Delta_v$) over $b\in \partial \Delta_w \backslash \Gamma'$ ($\phi(b)\in \partial \Delta_v \backslash \Gamma$). Then we have\\
\begin{th} 
For any $b\in \partial \Delta_w \backslash \Gamma'$, $Y_b$ is naturally dual to $X_{\phi(b)}$.\\
\end{th}
{\bf Proof:}
Based on above propositions, duality is very easy to establish. Only thing that need to be addressed is that duality defined in two ways according to $U_m$ or $\alpha_n^0$ for $b\in U_m\cap \alpha_n^0$ coincide. For this purpose, one only need to show that the following diagram commutes.\\
\[
\begin{array}{ccccccc}
X_{\phi(b)}&\rightarrow& T_m^\vee & \ \ &Y_b&\rightarrow& T_n^\vee \\
\downarrow&&\downarrow& \ &\downarrow&&\downarrow\\
T_n&\leftarrow& T_{\bf R}& \ &T_m&\leftarrow& T_{\bf R}^\vee\\
\end{array}
\]\\
This is proved in proposition \ref{ha} and \ref{hc}.
\begin{flushright} $\Box$ \end{flushright}
With explicit identification of fibres in place, monodromy computation becomes a piece of cake! Consider the path $\gamma_{nmn'm'} = \alpha_n^0 U_m \alpha_{n'}^0 U_{m'} \alpha_n^0$ on $\partial \Delta_v$, where $n,n'\in \partial \Delta^\vee$, $m,m'\in \Delta^0$ satisfying\\
\[
\langle m,n\rangle=\langle m,n'\rangle = \langle m',n\rangle = \langle m',n'\rangle = -1.
\]\\
This condition implies that $\alpha_n$ and $\alpha_{n'}$ have common face that contains $m,m'$. Correspondingly we have the diagram\\
\[
\begin{array}{ccc}
N_n&\rightarrow& N_m\\
\uparrow&&\downarrow\\
N_{m'}&\leftarrow& N_{n'}\\
\end{array}
\]\\
\[
x\rightarrow x+\langle m,x\rangle n \rightarrow x+\langle m,x\rangle n + \langle m',x+\langle m,x\rangle n\rangle n'= x+\langle m,x\rangle n + \langle m'-m,x\rangle n'
\]\\
Compose the four operators and modulo $n$, we get\\ 
\begin{th}
The monodromy operator along $\gamma_{nmn'm'}$ is\\
\[
[x] \rightarrow [x] + \langle m'-m,x\rangle [n']\ \ \ {\rm for}\ [x]\in N_n
\]\\
\end{th}
{\bf Remark:} Now we have find an extremely simple way to compute monodromy. All our monodromy computation in \cite{lag1} can be much easily performed by this method.\\\\
We can similarly compute monodromy along the corresponding path for the mirror Lagrangian torus fibration of the mirror Calabi-Yau hypersurface. Consider the corresponding path $\phi(\gamma_{nmn'm'}) = U_n \alpha_m^0 U_{n'} \alpha_{m'}^0 U_n$ on $\partial \Delta_w$. We have the corresponding diagram\\
\[
\begin{array}{ccc}
M_n&\rightarrow& M_m\\
\uparrow&&\downarrow\\
M_{m'}&\leftarrow& M_{n'}\\
\end{array}
\]\\
\[
y\rightarrow y+\langle y,n'\rangle m \rightarrow y+\langle y,n'\rangle m + \langle y+\langle y,n'\rangle m, n\rangle m'= y+\langle y,n'\rangle (m - m')
\]\\
Compose the four operators and modulo $n$, we get\\
\begin{th}
The monodromy operator along $\phi(\gamma_{nmn'm'})$ is\\
\[
y \rightarrow y+\langle y,n'\rangle (m - m')\ \ \ {\rm for}\ y\in M_n.
\]\\
\end{th}
Compare with the monodromy operator along $\gamma_{nmn'm'}$, we get\\
\begin{th}
The monodromy operator along $\phi(\gamma_{nmn'm'})$ is naturally dual to the monodromy operator along $\gamma_{nmn'm'}$.\\
\end{th}
{\bf Remark:} Although this duality result can be easily shown by the two monodromy formulas, more fundamentally, it is directly implied by the explicit description of fibres described in the previous propositions.\\\\ 
Summerize our results, we have proved the symplectic topological version of SYZ conjecture for quintic Calabi-Yau.\\
\begin{th}
For generic quintic Calabi-Yau $X$ near the large complex limit, and its mirror Calabi-Yau $Y$ near the large radius limit, there exist corresponding Lagrangian torus fibrations\\
\[
\begin{array}{ccccccc}
X_{\phi(b)}&\hookrightarrow& X& \ \ &Y_b&\hookrightarrow& Y\\
&&\downarrow& \ &&&\downarrow\\
&& \partial \Delta_v& \ &&& \partial \Delta_w\\
\end{array}
\]\\
with singular locus $\Gamma \subset \partial \Delta_v$ and $\Gamma' \subset \partial \Delta_w$, where $\phi:\partial \Delta_w \rightarrow \partial \Delta_v$ is a natural homeomorphism, $\phi(\Gamma')=\Gamma$ and $b\in \partial \Delta_w \backslash \Gamma'$. The corresponding fibres $X_{\phi(b)}$ and $Y_b$ are naturally dual to each other.\\
\end{th}
Our construction also give the dual relation of the singular fibres. Recall the discussion of singular fibres of generic Lagrangian torus fibration in \cite{lag3} section 7. For generic Lagrangian torus fibration, the singular fibres over smooth part of $\Gamma$ are of type $I$, singular fibres over verteces of $\Gamma$ are of type $II$ or type $III$.\\
\begin{th}
For $b\in \Gamma'$, if $Y_b$ is a singular fibre of type $I$, $II$, $III$, then $X_{\phi(b)}$ is a singular fibre of type $I$, $III$, $II$.\\
\end{th}
{\bf Remark:} The last piece of the SYZ puzzle we have not yet discussed is the construction of a section of the Lagrangian fibration. With the explicit description of Lagrangian fibres in this section, it is not hard to construct the section. Roughly, one can take the identity section on each piece with explicit description. When all the coefficients are real, they almost piece together to form an approximate global section, with error depending on how close the Calabi-Yau is near the Lagrangian complex limit. With a more careful construction, we will get a global section. In the case of general coefficients, similar construction will give us a canonical global section if understood in a suitable sense. We will describe the precise construction of gloal section in \cite{mon}, where more generally we will discuss monodromy near the large complex limit of general Calabi-Yau hypersurface in toric variety.\\\\

\se{Construction of mirror manifold via Lagrangian torus fibration}
As we mentioned in the introduction, the original SYZ mirror conjecture is rather sketchy in nature, with no mentioning of singular locus, singular fibres and duality of singular fibres, which is essential if one wants, for example, to use SYZ to construct mirror manifold. Our discussion of construction of Lagrangian torus fibration and symplectic topological SYZ of generic Calabi-Yau hypersurface in toric variety explicitly produce for us the 3 types of generic singular fibres (type $I$, $II$, $III$ as described in \cite{lag3}) and the way they are dual to each other under the mirror symmetry. This together with the knowledge of singular locus from our construction will enable us to give a more precise formulation of SYZ nirror conjecture. This precise formulation will naturally suggest a way to construct mirror manifold in general starting from a generic Lagrangian torus fibration of a Calabi-Yau manifold.\\\\
{\bf Precise SYZ mirror conjecture}
For any Calabi-Yau 3-fold $X$, with Calabi-Yau metric $\omega_g$ and holomorphic volume form $\Omega$, there exists a special Lagrangian fibration of $X$ over $S^3$\\
\[
\begin{array}{ccc}
T^3&\hookrightarrow& X\\
&&\downarrow\\
&& S^3
\end{array}
\]\\
with a special Lagrangian section and codimension 2 singular locus $\Gamma\subset S^3$, such that general fibres (over $S^3\backslash \Gamma$) are 3-torus. For generic such fibration, $\Gamma$ is a graph with only 3-valent verteces. Let $\Gamma = \Gamma^1\cup \Gamma^2 \cup \Gamma^3$, where $\Gamma^1$ is the smooth part of $\Gamma$, $\Gamma^2 \cup \Gamma^3$ is the set of the verteces of $\Gamma$. For any leg $\gamma\subset \Gamma^1$, the monodromy of $H_1(X_b)$ of fibre under suitable basis is\\
\[
T_\gamma =\left(
\begin{array}{ccc}
1 &  1 & 0\\
0 &  1 & 0\\
0 &  0 & 1
\end{array}
\right)
\]\\
Singular fibre along $\gamma$ is of type $I$.\\\\
Consider a vertex $P\in\Gamma^2 \cup \Gamma^3$ with legs $\gamma_1$, $\gamma_2$, $\gamma_3$. Correspondingly, we have monodromy operators $T_1$, $T_2$, $T_3$.\\\\
For $P\in\Gamma^2$, under suitable basis we have\\
\[
T_1 =\left(
\begin{array}{ccc}
1 &  1 & 0\\
0 &  1 & 0\\
0 &  0 & 1
\end{array}
\right)\ \ 
T_2 =\left(
\begin{array}{ccc}
1 &  0 & -1\\
0 &  1 & 0\\
0 &  0 & 1
\end{array}
\right)\ \ 
T_3 =\left(
\begin{array}{ccc}
1 &  -1 & 1\\
0 &  1 & 0\\
0 &  0 & 1
\end{array}
\right)
\]\\
Singular fibre over $P$ is of type $II$.\\\\  
For $P\in\Gamma^3$, under suitable basis we have\\
\[
T_1 =\left(
\begin{array}{ccc}
1 &  0 & 0\\
1 &  1 & 0\\
0 &  0 & 1
\end{array}
\right)\ \ 
T_2 =\left(
\begin{array}{ccc}
1 &  0 & 0\\
0 &  1 & 0\\
-1 &  0 & 1
\end{array}
\right)\ \ 
T_3 =\left(
\begin{array}{ccc}
1 &  0 & 0\\
-1 &  1 & 0\\
1 &  0 & 1
\end{array}
\right)
\]\\
Singular fibre over $P$ is of type $III$.\\\\
The special Lagrangian fibration for the mirror Calabi-Yau manifold $Y$ has the same base $S^3$ and singular locus $\Gamma \subset S^3$. For $b\in S^3\backslash \Gamma$, $Y_b$ is the dual torus of $X_b\cong T^3$. In another word, the $T^3$-fibrations\\
\[
\begin{array}{ccc}
T^3&\hookrightarrow& X\\
&&\downarrow\\
&& S^3\backslash\Gamma
\end{array}
\ \ \ \ \ \ 
\begin{array}{ccc}
T^3&\hookrightarrow& Y\\
&&\downarrow\\
&& S^3\backslash\Gamma
\end{array}
\]\\
are dual to each other. In particular the monodromy operator will be dual to each other.\\\\
For the fibration of $Y$, singular fibres over $\Gamma^1$ should be type $I$, singular fibres over $\Gamma^2$ should be type $III$, singular fibres over $\Gamma^2$ should be type $II$. Namely, dual singular fibre of a type $I$ singular fibre is still type $I$. Type $II$ and $III$ singular fibres are dual to each other.
\begin{flushright} $\Box$ \end{flushright}
Our work (in \cite{lag3,ci} and this paper) proved the symplectic topological version of this precise SYZ mirror conjecture for quintic, \cy hypersurface and more generally \cy complete intersection in toric varieties.\\\\
{\bf Remark:} Although it is widely believed that the special Lagangian fibration is a $C^\infty$-map and the singular locus is of codimension 2. It is still possible that the singular locus of the actual special Lagrangian is of codimension 1 (as the Lagrangian fibration we constructed before we modify the singular locus to codimension 2). If that is the case, the special Lagangian fibration will not be a $C^\infty$-map. Even if this is the case, we believe the symplectic topological version of this precise SYZ conjecture will still hold true, which is sufficient for us to construct the mirror manifold symplectic topologically.
\begin{flushright} $\Box$ \end{flushright}
This precise formulation of SYZ conjecture naturally suggests a way to construct mirror manifold $Y$ in general starting from a generic \l torus fibration of a \cy manifold $X$. The mirror manifold $Y$ can be constructed as the dual \l torus fibration over $S^3\backslash \Gamma$. Filling singular fibres is a delicate issue. Roughly speaking, dual singular fibre of a type $I$ singular fibre is still type $I$. Type $II$ and $III$ singular fibres are dual to each other.\\\\
To make the construction of mirror manifold precise, it is necessary to construct explicit local model for type $II$ and $III$ singular fibres. (Our gradient flow construction naturally provide this kind of local models.) Then one need to fit things together using the horizontal section and horizontal foliation determined by the Lagrangian fibration and the horizontal section. We intend to discuss this construction in greater detail in the future.\\\\

Current address:\\
Wei-Dong Ruan\\
Department of Mathematics\\
University of Illinois at Chicago\\
Chicago, IL 60607-7045\\\\
email: ruan@math.uic.edu


\begin{thebibliography}{11}
\bibitem{AGM}
Aspinwall, P.S., Greene, B.R., Morrison, D.R.,
``The Monomial-Divisor Mirror Map",
{\it Inter. Math. Res. Notices} 12 (1993), 319-337.
%
\bibitem{Can}
Candelas, P., de la Ossa, X.C., Green, P., Parkes, L.,
``A Pair of Calabi-Yau Manifolds as an Exactly Soluble Superconformal Theory",
in {\it Essays on Mirror Symmetry}, edited by S.-T. Yau.
%
\bibitem{Gross1}
Gross, M.,
``Special Lagrangian Fibration I: Topology",
alg-geom 9710006.
%
\bibitem{Gross2}
Gross, M.,
``Special Lagrangian Fibration II: Geometry",
alg-geom 9809072.
%
\bibitem{Gross3}
Gross, M.,
``Topological Mirror Symmetry",
alg-geom 9909015.
%
\bibitem{GW}
Gross, M. and Wilson, P.M.H.,
``Mirror Symmetry via 3-torus for a class of Calabi-Yau Threefolds",
to appear in {\it Math. Ann.}
%
\bibitem{HL}
Harvey, R. and Lawson, H.B.,
``Calibrated Geometries",
{\it Acta Math.} 148 (1982), 47-157.
%
\bibitem{H}
Hitchin, N.,
``The Moduli Space of Special Lagrangian Submanifolds",
dg-ga 9711002
%
\bibitem{lag1}
Ruan, W.-D.,
``Lagrangian torus fibration of quintic Calabi-Yau hypersurfaces I: Fermat type quintic case",
dg-ga 9904012.
%
\bibitem{lag2}
Ruan, W.-D.,
``Lagrangian torus fibration of quintic Calabi-Yau hypersurfaces II: Technical results on gradient flow construction",
(preprint).
%
\bibitem{lag3}
Ruan, W.-D.,
``Lagrangian torus fibration of quintic Calabi-Yau hypersurfaces III: Symplectic topological SYZ mirror construction for general quintics",
dg-ga 9909126.
%
\bibitem{N}
Ruan, W.-D.,
``Newton polygon, string diagram and toric variety",
(preliminary version).
%
\bibitem{ci}
Ruan, W.-D.,
``Lagrangian torus fibration and mirror symmetry of Calabi-Yau complete intersections in toric variety",
(preliminary version).
%
\bibitem{smooth}
Ruan, W.-D.,
``Smoothing of Lagrangian fibration map",
(In preparation).
%
\bibitem{mon}
Ruan, W.-D.,
``Monodromy near the large complex limit and horizontal sections of the Lagrangian torus fibration of Calabi-Yau hypersurface in toric variety",
(In progress).
%
\bibitem{SYZ}
Strominger, A.,Yau, S.-T. and Zaslow, E, 
``Mirror Symmetry is T-duality",
{\it Nuclear Physics} B 479 (1996),243-259.
%
\bibitem{Z}
Zharkov, I.,
``Torus Fibrations of Calabi-Yau Hypersurfaces in Toric Varieties and Mirror Symmetry",
alg-geom 9806091
\end{thebibliography}
\end{document}